\documentclass{article}

\usepackage{arxiv}
\usepackage{float}
\usepackage[utf8]{inputenc}
\usepackage[T1]{fontenc}
\usepackage{hyperref}
\usepackage{url}
\usepackage{booktabs}
\usepackage{subcaption}
\usepackage{amsfonts}
\usepackage{nicefrac}
\usepackage{microtype}
\usepackage{amsmath}
\usepackage{amssymb}
\usepackage{amsthm}
\usepackage{cleveref}
\usepackage{lipsum}
\usepackage{graphicx}
\usepackage[numbers]{natbib}
\usepackage{doi}
\usepackage{geometry}
\usepackage{xfrac}
\usepackage[group-separator={,}]{siunitx}

\theoremstyle{plain}
\newtheorem{theorem}{Theorem}[section]

\theoremstyle{definition}

\newtheorem{remark}[theorem]{Remark}

\title{Summation-by-Parts Finite-Difference Method for Linear Shallow Water Equations on Staggered Curvilinear Grids in Closed Domains}

\author{
    V. Shashkin, G. Goyman, I. Tretyak \\
G.I. Marchuk Institute of Numerical Mathematics, RAS, Moscow, Russia  \\  Moscow Center for Advanced Studies, Moscow, Russia \\ Hydrometeorological center of Russia, Moscow, Russia \\
v.shashkin@inm.ras.ru
}

\renewcommand{\shorttitle}{SBP FD for Linear Shallow Water Equations on Staggered Curvilinear Grids}

\begin{document}
\maketitle

\begin{abstract}
This work focuses on developing high-order energy-stable schemes for wave-dominated problems in closed domains using staggered finite-difference summation-by-parts (SBP FD) operators. We extend the previously presented uniform staggered grid SBP FD approach to non-orthogonal curvilinear multi-block grids and derive new higher-order approximations. The combination of Simultaneous-Approximation-Terms (SAT) and projection method is proposed for the treatment of interface conditions on a staggered grid. This reduces approximation stiffness and mitigates stationary wave modes of pure SAT approach. Also, energy-neutral discrete Coriolis terms operators are presented. The proposed approach is tested using the linearized shallow water equations on a rotating sphere, a testbed relevant for ocean and atmospheric dynamics. Numerical experiments show significant improvements in capturing wave dynamics compared to collocated SBP FD methods.
\end{abstract}

\keywords{
Summation-by-Parts \and Interface conditions \and Staggered grids \and Shallow water model
}

\section{Introduction}
Accurate simulations of wave dynamics are crucial for various applications, including ocean and atmosphere modeling, seismology, and electrodynamics. These fields often involve complex phenomena where wave interactions and dispersive effects play a crucial role, requiring numerical methods to effectively handle these intricate processes.

To accurately model wave propagation, an effective spatial approximation must ensure a high order of approximation, stability and an accurate representation of wave dispersion characteristics. This combination is crucial to avoid unphysical results and ensure accurate long-term simulations of wave dynamics.

A systematic approach to developing high-order energy stable schemes is based on the concept of summation-by-parts (SBP) operators \cite{kreiss1974finite,strand_sbp_1994,Olsson_projection_I,Olsson_projection_II,SBP_rev_2014,lundquist2024encapsulated}. 
The SBP approach aims to provide approximations of differential operators that satisfy a discrete version of the integration-by-parts formula. When coupled with a suitable boundary conditions treatment, the SBP method often results in semi-discrete approximations that possess an analogue of the energy conservation law, thereby ensuring the stability of numerical solutions.

In this work, we focus on finite-difference summation-by-parts (SBP FD) methods, which have undergone substantial advancements in recent years and are extensively utilized in various applications, including the modeling of wave-dominated problems \cite{SBP_rev_2014,svard2014review}. These methods enable the formulation of simple yet effective schemes on structured logically rectangular grids, enhancing the efficiency and accuracy of numerical simulations.

Standard SBP FD operators on uniform grids with collocated variables use ordinary central finite differences, with the exception of a limited number of points located near the boundaries, where one-sided differences are employed to preserve the SBP property \cite{kreiss1974finite,strand_sbp_1994,SBP_rev_2014}. Each SBP operator is associated with a norm matrix, whose structure is crucial as it affects both the approximation order and the stability of the developed schemes. For problems on curvilinear grids, diagonal norm matrices are generally preferred, as they are required for constructing provably stable schemes \cite{svard2004coordinate}. It is shown \cite{kreiss1974finite} that diagonal norm SBP operators with a $2s$-order accurate repeating interior stencil have a boundary closure order that does not exceed $s$. This limits the global $l_2$-norm error convergence rate to $s+1$ for systems of first-order hyperbolic equations \cite{gustafsson1975convergence}.

SBP FD methods can be naturally extended to multi-block geometries while preserving the key properties of single-block implementations, such as stability and high order of approximation. Utilizing a multi-block mesh facilitates handling of complex geometries, enables the delineation of areas with different physical properties, and provides opportunities for static mesh adaptation, where local mesh resolution can be increased in regions of interest without incurring unnecessary computational costs in less critical areas \cite{mattsson2010stableMultiblock,almquist2019order,tretyak2023multiresolution}. 

To ensure stability when employing SBP FD, particular diligence is required when imposing boundary conditions (and interface conditions in the case of multi-block grids). A widely adopted approach to address these challenges involves the combination of SBP operators with simultaneous approximation terms (SBP-SAT) \cite{carpenter1994time}. This methodology introduces penalty terms to the right-hand sides of the equations, enabling weak enforcement of boundary and interface conditions while preserving discrete analogues of conservation laws.  

Another strategy for managing interface and boundary conditions is the projection method (SBP-P) \cite{Olsson_projection_I,Olsson_projection_II,mattsson2018improvedProjection}. This approach employs an orthogonal projection to restrict the numerical solution to a subspace of grid functions that exactly satisfies the boundary or interface conditions. When combined with SBP operators, this approach yields energy-conservative approximations. While both SBP-SAT and SBP-P methods enable the construction of stable high-order schemes, the characteristics of the resulting approximations can vary significantly. In particular, issues such as the stiffness of the approximation and the emergence of spurious computational modes may occur. 

In the context of first-order hyperbolic equations systems, high-order schemes with collocated variables (unstaggered grids) may not adequately capture wave propagation phenomena across all spatial scales. Specifically, the use of central difference approximations can lead to phase velocity attenuation for the shortest wavelength components of the solution, as well as spurious inversions of group velocity. In contrast, methods using grids with variables staggering (staggered grids) mitigate these issues, generally resulting in more accurate solutions for wave-dominated problems \cite{ArakawaLamb1977,Randall_GeostAdj,Thuburn2011_hor_basic_ideas,zingg2000comparison,GoymanShashkinJCP2021}. This observation has motivated the development of staggered SBP FD schemes.

Diagonal norm staggered summation-by-parts finite difference operators were first introduced in \cite{OReilly_staggered_sbp_2017}, where they were combined with SAT boundary treatment to solve 2D wave equations. In subsequent work \cite{OReilly_staggered_sbp_2020}, these operators, along with specialized interpolation procedures, were used to solve the acoustic wave equation on staggered curvilinear grids. 
A distinctive feature of the proposed operators is the use of collocated nodes at domain boundaries, allowing to utilize the same form of SBP property as in the unstaggered grid case. However, this approach makes points placement for one of variables non-uniform, with two additional nodes shifted by a half grid spacing.
The use of non-uniform grid makes the approach \cite{OReilly_staggered_sbp_2017} less attractive for several reason: first, the use of extra grid nodes results in disparity between vector and scalar degrees of freedom, which may lead to the emergence of spurious computational branches in the wave solutions \cite{walters1983analysis,le2005some,Stan_grid_rev}. Secondly, uniform grids generally lead to simpler, more efficient and faster algorithms, which is important when using this approach in real-world models. 

The work \cite{Gao_staggered_sbp_2019} considers staggered SBP FD operators of order 4/2 (interior/near-boundary) utilizing conventional uniform grids. The proposed approximations combined with the SAT approach were used to solve 2d acoustic equations on a multi-block grid with variable resolution. However, only Cartesian grids were considered. The extension to non-orthogonal curvilinear meshes requires the derivation of interpolation operators between locations of staggered variables that are compatible with the norm matrices of SBP FD operators. 
Additionally, it necessitates the construction of a metric tensor approximation that preserves both symmetry and positive definiteness, as outlined in \cite{OReilly_staggered_sbp_2020}. Furthermore, in \cite{Gao_staggered_sbp_2019} it was noticed that the use of the SBP-SAT approach leads to a noticeably more stringent limitation on the time step than expected, which also requires further research and possibly modification of the interface conditions imposing procedure.

Overall, the primary goal of this work is to address the problem of unsatisfactory wave dispersion properties of collocated grid SBP finite-difference schemes at small spatial scales. The way to solve this problem is improving the staggered SBP finite-difference methodology first derived in works \cite{OReilly_staggered_sbp_2017}, \cite{OReilly_staggered_sbp_2020}, \cite{Gao_staggered_sbp_2019}. As compared to previous works we:

- extend the energy-stable uniform grid staggered SBP scheme \cite{Gao_staggered_sbp_2019} to multi-block curvilinear grid case;

- formulate the hybrid SAT-projection procedure to impose interface conditions at grid block boundaries;

- derive sixth order accurate staggered SBP-FD  approximations for $d/dx$;

- derive energy-conserving Coriolis terms approximation consistent with staggered SBP-discretization.

The first two points contribute to the reproduction of the wave-dispersion properties by removing spurious fast and stationary modes associated with grid-blocks interface treatment and also remove excessive problem stiffness noted in \cite{Gao_staggered_sbp_2019}. The third point improves accuracy of staggered SBP-FD. Fourth point is crucial for application of staggered SBP-FD in the geophysical hydrodynamics field.

For the case studies, we focus on the numerical solution of the linearized shallow water equations on the rotating sphere 

\begin{subequations}\label{eq:shallow_water}
    \begin{equation}
    \frac{\partial \vec v}{\partial t} = f \vec v^\perp -g\nabla h,
    \end{equation}
    \begin{equation}
        \frac{\partial h}{\partial t} = -\mathcal{H} \nabla\cdot\vec v,
    \end{equation}
\end{subequations}
where $\vec v$ is the flow velocity vector, $f$ is the Coriolis parameter, $\vec v^\perp = k\times \vec v$, where $k$ is the vertical unit vector, $g$ is the gravity acceleration, $h$ is the fluid height perturbation, $\mathcal{H}$ is the fluid mean height. Note that system of equations (\ref{eq:shallow_water}) with zero $f$ are analogous to 2D acoustic/elastic waves equations and some Maxwell equations approximations.

The linearized shallow water equations serve as an important testbed for evaluating new numerical methods intended for use in ocean and atmospheric models \cite{Williamson_tests,Galewsky_test}. These equations offer a simplified yet representative framework for the isolated simulation of inertia-gravity and Rossby waves dynamics, allowing for the assessment of numerical solution characteristics that are critical for accurately capturing the evolution of the ocean and atmosphere. Key characteristics of interest among others include the accuracy in reproducing the phase and group velocities, the absence of non-physical computational modes, and the reproduction of geostrophic adjustment.

The article is organized as follows: Section \ref{sec:sbp_sat_1d} presents the staggered SBP finite difference (FD) operators in one dimension. In Section \ref{sec:2d_sbp_sat}, we describe the 2D staggered SBP FD approximation on multi-block grids. Section \ref{sec:model formulation} outlines the formulation of the test problem and its approximation, while Section \ref{sect:mimetic_prove} provides proofs of the mimetic properties associated with the considered SBP FD approximation. Section \ref{sect:num_exp} analyzes the results of numerical experiments, and Section \ref{sec:conclusion} concludes with a summary of the findings.
\section{Introduction to staggered SBP FD method in one dimension}\label{sec:sbp_sat_1d}
\subsection{Staggered SBP FD first derivative operators}
We begin by describing the 1D staggered SBP FD methods as presented in \cite{Gao_staggered_sbp_2019}.
Consider two uniform grids defined on the domain $[a,b]$
\begin{align}\label{eq:xc_xv}
& \mathbf{x}^v=(x^v_1,\ldots,x^v_{N+1})^T, &  x^v_{i}&=a+(i-1)\Delta x, &i&\in[1,N+1], \\
& {\mathbf{x}}^c=(x^c_1,\ldots,x^c_N)^T,   &  x^c_i &= a + \left(i-\frac{1}{2}\right)\Delta x, &i&\in[1,N], 
\end{align}
with $\Delta x = (b-a)/N$. Grid points $x_i^c$ correspond to cell centers, and $x_i^v$ correspond to cell vertices. Then, column vectors
 \begin{equation*}
 {\mathbf{h}}=({h}_1, \ldots,{h}_{N+1})^{\mathrm{T}}, \quad \mathbf{u}=(u_1, \ldots,u_{N})^{\mathrm{T}}
 \end{equation*}
represent the values of some functions $h(x)$ and $u(x)$ at the grids $\mathbf{x}^v$ and ${\mathbf{x}^c}$, respectively.
 
We define matrices $D_{vc}$ and $D_{cv}$ approximating the first derivative. The matrix $D_{vc}$ approximates the derivative at $\mathbf{x}^c$ points using function values at ${\mathbf{x}}^v$ points, and vice versa for the $D_{cv}$ matrix:
 \begin{equation*}
     \frac{\partial u}{\partial x}\approx D_{cv} \mathbf{u}, \quad 
     \frac{\partial h}{\partial x}\approx D_{vc} {\mathbf{h}}.
 \end{equation*}

The SBP property is a discrete analogue of continuous integration by parts:
\begin{equation}
    \int_{a}^b u\frac{\partial h}{\partial x} dx = 
    -\int_{a}^b h\frac{\partial u}{\partial x} dx + u(b)h(b)-u(a)h(a). 
\end{equation}

Suppose $H_c$ and $H_v$ are diagonal positive definite matrices that define discrete quadrature rules on the $\mathbf{x}^c$ and $\mathbf{x}^v$ grids, respectively. Then, the SBP property is written as:
\begin{equation}\label{eq:sbp_property}
    \mathbf{u}^T H_c D_{vc} {\mathbf{h}} = -\mathbf{{h}}^T {H}_v D_{cv} \mathbf{u} +
    h_{N+1} \mathbf{r}^T {\mathbf{u}} - h_1 \mathbf{l}^T {\mathbf{u}},
\end{equation}
where $\mathbf{l}^T$ and $\mathbf{r}^T$ are row vectors such that $\mathbf{l}^T{\mathbf{u}}\approx u(a)$ and $\mathbf{r}^T{\mathbf{u}}\approx u(b)$ are extrapolations to the left and right domain boundaries, respectively. Extrapolation is necessary because the ${\mathbf{x}}^c$ grid does not contain boundary points. It is convenient to use the following notation for the boundary-flux terms:
\begin{equation}
     h_{N+1} \mathbf{r}^T {\mathbf{u}} -   h_1 \mathbf{l}^T {\mathbf{u}} \equiv \mathbf{h}^T (\mathbf{e}_r\mathbf{r}^T- \mathbf{e}_l\mathbf{l}^T) {\mathbf{u}}\equiv \mathbf{h}^T R_{cv} {\mathbf{u}},
\end{equation}
where $\mathbf{e}_r = (1,0,...,0)^T \in \mathbb{R}^{(N+1)\times1}$, $\mathbf{e}_l = (0,...,0,1)^T \in \mathbb{R}^{(N+1)\times1}$. 

A pair of differentiation matrices $D_{cv}$, $D_{vc}$ that satisfies the condition (\ref{eq:sbp_property}) for some $H_c$, $H_v$, $R_{cv}$ matrices, will be referred to as staggered SBP FD operators.  

The SBP FD operators coincide with conventional central finite differences at grid points away from boundaries
\begin{equation}
    \frac{\mathrm{d} h}{\mathrm{d} x} \bigg |_{x=x^v_i} = \frac{u_{i} - u_{i-1}}{\Delta x} + O(\Delta x^2),
\end{equation}
\begin{equation}
    \frac{\mathrm{d} h}{\mathrm{d} x} \bigg |_{x=x^v_i} = \frac{27(u_{i} - u_{i-1})-(u_{i+1} - u_{i-2})}{24\Delta x} + O(\Delta x^4),
\end{equation}
while one-sided or extrapolation-like difference formulae are used near the boundaries.
The use of diagonal quadrature matrices $H_c$ and $H_v$ implies that $2s$-order accurate operators at the inner points can simultaneously have only $s$-order accuracy near the edges. In this case, the boundary extrapolation operators $\mathbf{l}$ and $\mathbf{r}$ need to have the $(s+1)$-order of accuracy.

In this work, we use $2/1$, $4/2$, and $6/3$  (interior/near the boundaries) order-accurate staggered SBP FD operators. 
The $2/1$ operators are defined as
\begin{equation}
	\mathrm H_v = \Delta x
		\begin{bmatrix}
			\sfrac{1}{2} &   &       &   &   \\
			& 1 &       &   &   \\
			&   &\ddots &   &   \\
			&   &       & 1 &   \\
			&   &       &   & \sfrac{1}{2} 
		\end{bmatrix}
	, \quad
 	\mathrm D_{cv} = \frac{1}{\Delta x}
		\begin{bmatrix}
			-1  & 1   &       &     \\
			-1 & 1 &       &     \\
			&\ddots   &\ddots &     \\
			&   &       -1& 1   \\
			&   &       -1&  1  
		\end{bmatrix},
\end{equation}

\begin{equation}
\mathrm H_c = \Delta x
		\begin{bmatrix}
			1 &   &       &   &   \\
			&   &\ddots &   &   \\
			&   &       &   & 1 
		\end{bmatrix}	
	, \quad
		\mathrm D_{vc} = \frac{1}{\Delta x}
		\begin{bmatrix}
			-1  & 1   &       &    \\
			      & \ddots & \ddots       &    \\
			&   &     -1  & 1   
		\end{bmatrix}		,
\end{equation}
\begin{equation}
\mathbf{l}^T =(\sfrac{3}{2}, -\sfrac{1}{2}, 0, \dots, 0)^T 	
	, \quad
		\mathbf{r}^T =(0, \dots,0,-\sfrac{1}{2}, \sfrac{3}{2})^T.
\end{equation}

The $4/2$ operators were first presented in \cite{Gao_staggered_sbp_2019} and for completeness of presentation are also given in the \ref{app:coefficients:quad:42}, \ref{app:coefficients:diff:42}.

To our knowledge, $6/3$-order staggered SBP FD operators have never been presented in the literature, so we derived them using the procedure described in \cite{SBP_rev_2014,Gao_staggered_sbp_2019}. 

First, we fix fourth-order boundary extrapolation operators:
\begin{equation}
    \mathbf{l} = \frac{1}{16}(35, -35, 21, -5,0,...,0)^T, \quad \mathbf{r} = \frac{1}{16}(0,...,0,-5,21,-35,35)^T.
\end{equation}
Imposing the approximation order and SBP property, we obtain the first derivative operators with two free parameters. These parameters can be used to optimize the operators in various ways, such as minimizing the spectral radius or truncation error at near-boundary points \cite{SBP_rev_2014,Diener_SBP_cubsph_like}. 

In this work, we consider two methods to choose these parameters. The first method follows the conventional approach of minimizing polynomial differentiation error. Specifically, the free parameters are determined by minimizing the objective function:
\begin{equation}\label{eq:poly_opt_63}
\left(\mathbf{Err}^4_{D_{cv}}\right)^T \mathbf{Err}^4_{D_{cv}} + \left(\mathbf{Err}^4_{D_{vc}}\right)^T \mathbf{Err}^4_{D_{vc}},
\end{equation}
where
\begin{align}
\mathbf{Err}^4_{D_{cv}} &= D_{cv} \left((x_1^c)^4, \dots, (x_N^c)^4 \right)^T -4\left((x_1^v)^3,\dots, (x_{N+1}^v)^3 \right)^T, \\
\mathbf{Err}^4_{D_{vc}} &= D_{vc} \left((x_1^v)^4, \dots, (x_{N+1}^v)^4 \right)^T-4\left((x_1^c)^3,\dots, (x_N^c)^3 \right)^T
\end{align}
are column vectors representing fourth-order polynomials differentiation error. According to \cite{Diener_SBP_cubsph_like}, minimizing operator's truncation error is closely related to minimizing the spectral radius. 

The second option is based on minimizing the error in calculating the second derivative for waves with wavelengths of $4\Delta x$ and $8 \Delta x$. The objective function is defined as
\begin{equation}\label{eq:wave_opt_63}
\overline{\mathbf{Err}^{4\Delta x}}^T\mathbf{Err}^{4\Delta x} + \overline{\mathbf{Err}^{8\Delta x}}^T\mathbf{Err}^{8\Delta x},
\end{equation}
where the error vectors are given by:
\begin{equation}
\begin{gathered}
    \mathbf{Err}^{k\Delta x} = (k\Delta x)^2 D_{cv}D_{vc}\mathbf{t}^{k\Delta x} - \mathbf{t}^{k\Delta x}, \\
    (\mathbf{t}^{k\Delta x})_m = \exp\left(i \frac{2\pi}{k\Delta x} x^v_m\right),
\end{gathered}
\end{equation}
and $\overline{(\ )}$ denotes complex conjugate operation. This optimization approach aims to reduce errors related to non-physical wave reflections at grid blocks interfaces. The selection of this objective function is empirical, with the goal of minimizing distortions near boundaries in the eigenvector structure of the discrete Laplace operator across both the initial and central regions of the wavenumber range. See Section \ref{sect:num_exp:poorly_resolved} for more details.

Coefficients of the $6/3$-order operators and the values of free parameters used in this work are given in the \ref{app:coefficients:quad:63}, \ref{app:coefficients:diff:63}.

\begin{remark}
        In derivation the differentiation operators here and in \cite{Gao_staggered_sbp_2019}, pre-defined extrapolation operators with a minimal stencil are used. However, the boundary stencil widths of the resulting differentiation operators generally allow for using wider extrapolation stencils. This enables the introduction of free parameters into the extrapolation operators, which can be used for further optimization of the differentiation matrices. This can be investigated in future works.
\end{remark}

We will also need interpolation operators between the $\mathbf{x}^v$ and $\mathbf{x}^c$ grids to transform covariant to contravariant components of the wind velocity vector and to approximate Coriolis force terms in the momentum equation. These operators have $2s$-order accuracy in the interior of the domain and $s$-order accuracy at the boundary. We introduce $P_{vc}$ and $P_{cv}$ operators, which interpolate grid function values from the $\mathbf{x}^v$-grid to the $\mathbf{x}^c$-grid and vice versa. It will be demonstrated below that these interpolation operators must meet the so-called SBP-preserving or inner product preserving property \cite{OReilly_staggered_sbp_2020}:
\begin{equation}\label{eq:sbp_interp_requirement}
    \mathbf{u}^T H_c P_{vc} {\mathbf{h}} = \mathbf{{h}}^T {H}_v P_{cv} \mathbf{u},
\end{equation}
to ensure the stability and mimetic features of the proposed discretizations. This property indicates that when using the SBP quadrature to evaluate an integral, the result should be the same regardless of whether $\mathbf{u}$ or $\mathbf{h}$ grid function values are interpolated. Note that SBP property preserving interpolation operators are constructed using the same quadrature matrices $H_c$ and $H_v$ that correspond to SBP FD operators $D_{cv}$ and $D_{vc}$.

The 2/1-order SBP-preserving interpolation operators are defined uniquely as follows:
\begin{equation}
	\mathrm P_{vc} = \frac{1}{2}
		\begin{bmatrix}
		  1	& 1     &     &    &   \\
		   & 1     & 1   &    &   \\
		   &       &  \ddots   & \ddots & \\
		   &       &          &  1 &  1 \\
		\end{bmatrix}
	, \quad
 	\mathrm P_{cv} = \frac{1}{2}
		\begin{bmatrix}
			2  & 0   &       &     \\
			1  & 1 &       &     \\
			      & \ddots   &\ddots &     \\
			&   &       1 & 1   \\
			&   &         &  2  
		\end{bmatrix},
\end{equation}.

After enforcing the SBP-preserving and approximation order conditions, the 4/2 and 6/3-order operators still have 2 and 6 free parameters, respectively. We select these parameters by minimizing the truncation errors of the operators. For the 4/2-order operator, we minimize
\begin{equation}\label{eq:interp_42_obj_func}
 \left(\mathbf{Err}^2_{P_{cv}}\right)^T \mathbf{Err}^2_{P_{cv}}+\left(\mathbf{Err}^2_{P_{vc}}\right)^T \mathbf{Err}^2_{P_{vc}}.   
\end{equation}
For the 6/3-order interpolation, we first minimize 
\begin{equation}\label{eq:interp_63_obj_func_1}
    \left(\mathbf{Err}^3_{P_{cv}}\right)^T \mathbf{Err}^3_{P_{cv}}+\left(\mathbf{Err}^3_{P_{vc}}\right)^T \mathbf{Err}^3_{P_{vc}}
\end{equation}
 and then minimize 
 \begin{equation}\label{eq:interp_63_obj_func_2}
     \left(\mathbf{Err}^4_{P_{cv}}\right)^T \mathbf{Err}^4_{P_{cv}}+\left(\mathbf{Err}^4_{P_{vc}}\right)^T \mathbf{Err}^4_{P_{vc}}.
 \end{equation}
Here, the interpolation error terms are defined as:
\begin{align}
    \mathbf{Err}^k_{P_{cv}} = P_{cv}\left( (x^c_1)^k,...,(x^c_N)^k \right)^T - \left( (x^v_1)^k,...,(x^v_{N+1})^k \right)^T, \\
    \mathbf{Err}^k_{P_{vc}} = P_{vc}\left( (x^v_1)^k,...,(x^v_{N+1})^k \right)^T - \left( (x^c_1)^k,...,(x^c_{N})^k \right)^T
\end{align}

The coefficients of the interpolation matrices and the values of the free parameters used in this work are provided in \ref{app:coefficients:interp}.

\begin{remark}
In Section \ref{sect:mimetic:Q}, we will demonstrate that a sufficient condition for the positive definiteness of the metric tensor is closely related to the spectral radius of the matrix \(P_{cv} P_{vc}\) and the skewness of the grid cells. Specifically, a higher degree of skewness in the grid cells necessitates that the spectral radius approaches unity. Therefore, in situations with highly skewed grids, an alternative strategy for parameter selection could entail minimizing the spectral radius of the matrix \(P_{cv} P_{vc}\). For instance, this minimization can be indirectly achieved by minimizing the Frobenius norm of the matrix.
\end{remark} 

\subsection{Imposing interface conditions}\label{subsec:boundary_conds_1d}

The summation-by-parts property is an essential component for correctly imposing various types of boundary conditions. When approximating equations on a multi-block mesh, we are particularly interested in constructing boundary conditions at the interfaces of the grid's blocks. To demonstrate how the SBP FD framework allows to do it, we first consider the 1D periodic domain $[a,\ b]$, so that $u(a)=u(b)$, $h(a)=h(b)$ and $x=a$, $x=b$ can be thought as an interface between two grid blocks.

First, note that application of SBP FD derivative operators results in different values at the points $x^v_1=a$, $x^v_{N+1}=b$ due to using one-sided stencils near the domain boundaries. Moreover, in general case $\mathbf{r}^T\mathbf{u} \neq \mathbf{l}^T\mathbf{u}$. All of this implies that the equivalent of the continuous relation
\begin{equation}\label{eq:IBP_periodic}
    \int_a^b u(x) \frac{\mathrm{d} h(x)}{dx} dx = -\int_a^b h(x) \frac{\mathrm{d} u(x)}{dx} dx
\end{equation}
does not hold for the discrete derivative operators $D_{cv}$, $D_{vc}$. To solve this problem, SBP FD operators can be combined with the simultaneous approximation terms (SAT) correction method \cite{carpenter1994time} or with the projection method \cite{Olsson_projection_I,Olsson_projection_II,mattsson2018improvedProjection}. 

\subsubsection{Simultaneous approximation terms (SAT) method}
The SAT technique involves modifying the derivative operators at near-boundary grid points in order to compensate for the terms in the summation-by-parts formula associated with the flux through the boundary. The SAT-corrected operators ${D}^S_{vc}$, ${D}^S_{cv}$ can be written as 
\begin{equation}\label{eq:1d_sat_vc}
    {D}^{S}_{vc} = D_{vc} -\frac{1}{2}H^{-1}_c(\mathbf{r} + \mathbf{l})(\mathbf{e}_r - \mathbf{e}_l)^T, 
\end{equation}
\begin{equation}\label{eq:1d_sat_cv}
    {D}^S_{cv} = D_{cv} -\frac{1}{2}H^{-1}_v(\mathbf{e}_r + \mathbf{e}_l)(\mathbf{r} - \mathbf{l})^T. 
\end{equation}
The penalty terms in SAT-correction are proportional to $h_{N+1}-h_1$ and $\mathbf{r}^T\mathbf{u} - \mathbf{l}^T\mathbf{u}$, and may be viewed as a weak enforcement of the interface condition. For the SAT-corrected operators, it can be shown that
\begin{equation}
\mathbf{u}^T H_c {D}^S_{vc} {\mathbf{h}} = -\mathbf{{h}}^T {H}_v {D}^S_{cv} \mathbf{u}.
\end{equation}

Although the SAT technique has been utilized effectively in a variety of works, we believe it has several drawbacks that are generally overlooked.

First, implementing interface conditions using the SAT technique can significantly increase the stiffness of difference operators compared to strong enforcement \cite{mattsson2018improvedProjection,eriksson2023boundary,eriksson2023non}. Specifically, this impact was seen in the context of staggered SBP FD as noted in \cite{Gao_staggered_sbp_2019}.
This increase in stiffness imposes more stringent restrictions on the maximum allowable time step in the case of explicit time integration, and slowdown of iterative solvers convergence for implicit methods.

The second issue is that there are additional degree of freedom (DoF) at the interface points when using the SAT technique. Thus, the use of staggered grid causes an imbalance in the DoFs for $\mathbf{h}$ and $\mathbf{u}$. Such an imbalance results in an extra zero eigenvalue for the discrete Laplace operator ($L=D^S_{cv}D^S_{vc}$) at the h-points, due to null space of the $D^S_{vc}$ matrix. For processes whose dynamics are governed by the Laplace operator (gravity waves, diffusion, etc.), this results in the formation of stationary computational modes that do not correspond to physical solutions \cite{walters1983analysis,le2005some}. This spurious mode can have an impact on the accuracy and stability of numerical solutions because its amplitude can grow as a result of a non-zero projection of floating-point arithmetic errors, nonlinear interactions, or external forcing in the right hand side of equations. Constructing an effective procedure for filtering such modes could be challenging, since artificial dissipation is usually based on the use of the same Laplace operator.

\subsubsection{SAT-projection method}\label{sect:sat_proj}
To avoid the problems described above, we propose to use an alternative approach for imposing interface conditions that combines projection and SAT approaches. To describe this method, let us first examine how the SAT correction of $D_{cv}$ modifies the summation-by-parts expression:
\begin{gather}
     \mathbf{u}^T H_c D_{vc} {\mathbf{h}} + \mathbf{{h}}^T {H}_v D^S_{cv} \mathbf{u} = \\
     \mathbf{u}^T H_c D_{vc} {\mathbf{h}} + \mathbf{{h}}^T {H}_v D_{cv} \mathbf{u} - \frac{1}{2} \mathbf{{h}}^T(\mathbf{e}_r + \mathbf{e}_l)(\mathbf{r} - \mathbf{l})^T \mathbf{{u}}
     = \\
     \frac{1}{2} \mathbf{{h}}^T(\mathbf{e}_r - \mathbf{e}_l)(\mathbf{r} + \mathbf{l})^T \mathbf{{u}} = (h_{N+1}-h_{1})\frac{(\mathbf{r} + \mathbf{l})^T }{2} \mathbf{u}
\end{gather}
The $D_{vc}$ and $D^S_{cv}$ matrices may be read as a pair of SBP FD operators, with the left and right boundary extrapolation operators equal to $\frac{(\mathbf{r} + \mathbf{l})^T }{2}$. Because boundary extrapolated values are the same, the interface term in the summation-by-parts formula vanishes if there is no jump of the $\mathbf{h}$ field values at the interface point. To achieve this condition, the projection method can be employed, which entails $l_2$-projecting the results of all operator computations at the grid $\mathbf{x}^v$ nodes onto the space of grid-functions that are continuous at the interface \cite{Olsson_projection_I,Olsson_projection_II}. Thus, the $D^S_{cv}$ and $D_{vc}$ are additionally modified as
\begin{equation}
    D^{P}_{cv} = AD_{cv}^S, \ D^{P}_{vc} = D_{vc} A.
\end{equation}
where $A$ is the projection matrix defined as \cite{Olsson_projection_I,Olsson_projection_II}:
\begin{equation}
    (A h)_i = 
    \begin{cases}
        \frac{(H_v h)_1 + (H_v h)_{N+1}}{(H_v)_{1,1} + (H_v)_{N+1,N+1}}, \ i=1 \ \text{or} \ i=N+1, \\
        h_i,\ \text{otherwise}.
    \end{cases}
\end{equation}
For the pair of matrices $D_{cv}^P$, $D_{vc}^P$ it can be shown
\begin{equation}
\begin{gathered}
     \mathbf{u}^T H_c D^P_{vc} {\mathbf{h}} + (A\mathbf{h})^T {H}_v D^P_{cv} \mathbf{u} = \\
     \mathbf{u}^T H_c D_{vc} {A\mathbf{h}} +
     \mathbf{(Ah)}^T {H}_v A D_{cv} \mathbf{u} 
     - \frac{1}{2} \mathbf{(Ah)}^T H_v A H_v^{-1}(\mathbf{e}_r + \mathbf{e}_l)(\mathbf{r} - \mathbf{l})^T \mathbf{{u}}
     = \\
     \mathbf{u}^T \left( H_c D_{vc} + D_{cv}^T H_v - \frac{1}{2}(\mathbf{r}-\mathbf{l})(\mathbf{e}_r+\mathbf{e}_l)^T\right)A\mathbf{h} = \\
     \frac{1}{2}\mathbf{u}^T (\mathbf{r}+\mathbf{l})(\mathbf{e}_r-\mathbf{e}_l)^TA\mathbf{h} = \frac{1}{2}\mathbf{u}^T (\mathbf{r}+\mathbf{l})\left((A\mathbf{h})_{N+1}-(A\mathbf{h})_{1}\right)\equiv 0,
     \end{gathered}
\end{equation}
where we have used orthogonal projection matrix properties $(H_v A)^T = H_v A$, $A^2 = A$. Therefore, $D_{cv}^P$ and $D_{vc}^P$ are SBP FD operators which mimic integration by parts formula for closed domains (\ref{eq:IBP_periodic}). Note that a hybrid SAT-projection approach on collocated grids was used in \cite{eriksson2023boundary,eriksson2023non} within a context of higher order hyperbolic equations, where the authors showed that the use of projection method leads to the smallest spectral radius of the spatial operators.

\subsubsection{Laplace spectrum analysis}\label{sect:laplace_spectrum}

\begin{figure}[h!]
     \centering
     \includegraphics[width=\textwidth]{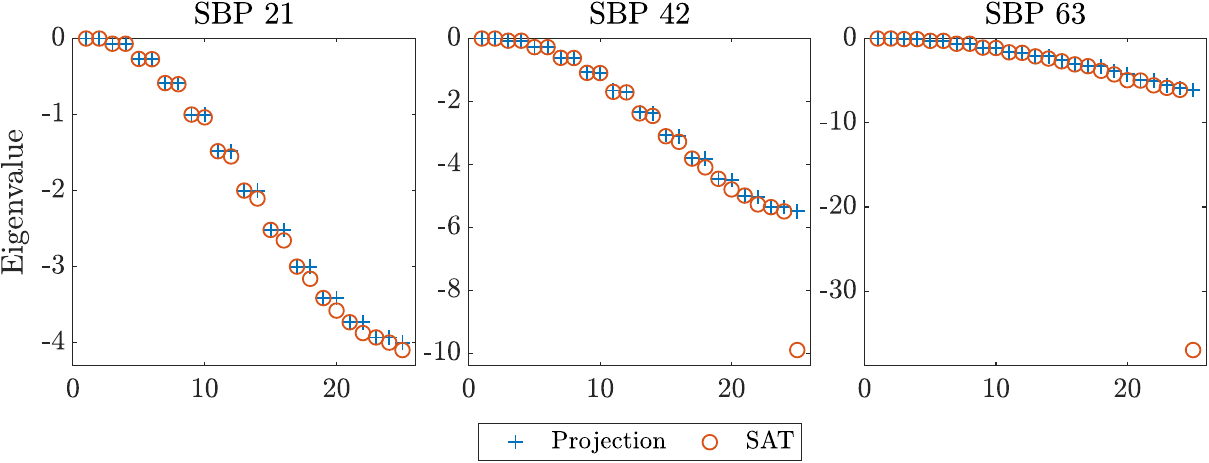}
     \caption{Eigenvalues of the discrete Laplace scaled by $\Delta x^2$ using SAT (red circles) or SAT-projection method (blue pluses) to impose interface conditions.}
     \label{fig:lap_spectrum}
\end{figure}

To show the difference of wave dynamics between SAT and SAT-projection SBP approximations we examine the spectrum of discrete Laplace operator. We consider matrices $L^S= D^S_{cv}D^S_{vc}$ and $L^P= D^P_{cv}D^P_{vc}$ for the 2/1, 4/2, and 6/3-order SBP FDs. Figure \ref{fig:lap_spectrum} shows the result of direct eigenvalues computation for the $L^S$ and $L^P$ matrices on the domain $[0,1]$ with $N=24$. In the case of SAT-based approximation, an outlier minimum eigenvalue exists for the 4/2 and 6/3-order approximations. The magnitude of the overshoot grows as the approximation order increases. These extreme eigenvalues correspond to eigenvectors with discontinuity at the interface points. This phenomenon is not observed when using the hybrid SAT-projection method, which appears to be explained by the fact that the projection operator filters out grid functions with discontinuities at the interface boundary. We also see that for the SAT-based approximation, two zero eigenvalues are observed, one of them corresponds to constant functions, and the second is spurious non-constant mode. For the hybrid SAT-projection method, an additional eigenvalue is also observed, but in this case it corresponds to the $(1,0,...,0,-1)^T$ eigenvector, and its amplitude is always zero for the interface-continuous grid functions, which is forced by the use of the projection operator and consistent initial conditions. 

Overall, the hybrid SAT-projection method effectively addresses the issues related to spurious fast and stationary modes that arise in the SAT-based approximation by filtering out grid functions with interface discontinuities. This approach ensures that only physically meaningful solutions are retained, enhancing the accuracy and stability of the numerical scheme.

\section{Two-dimensional multi-block SBP operators}\label{sec:2d_sbp_sat}

\begin{figure}
    \centering
    \includegraphics[trim={0 5cm 0 5cm},clip,width=0.7\linewidth]{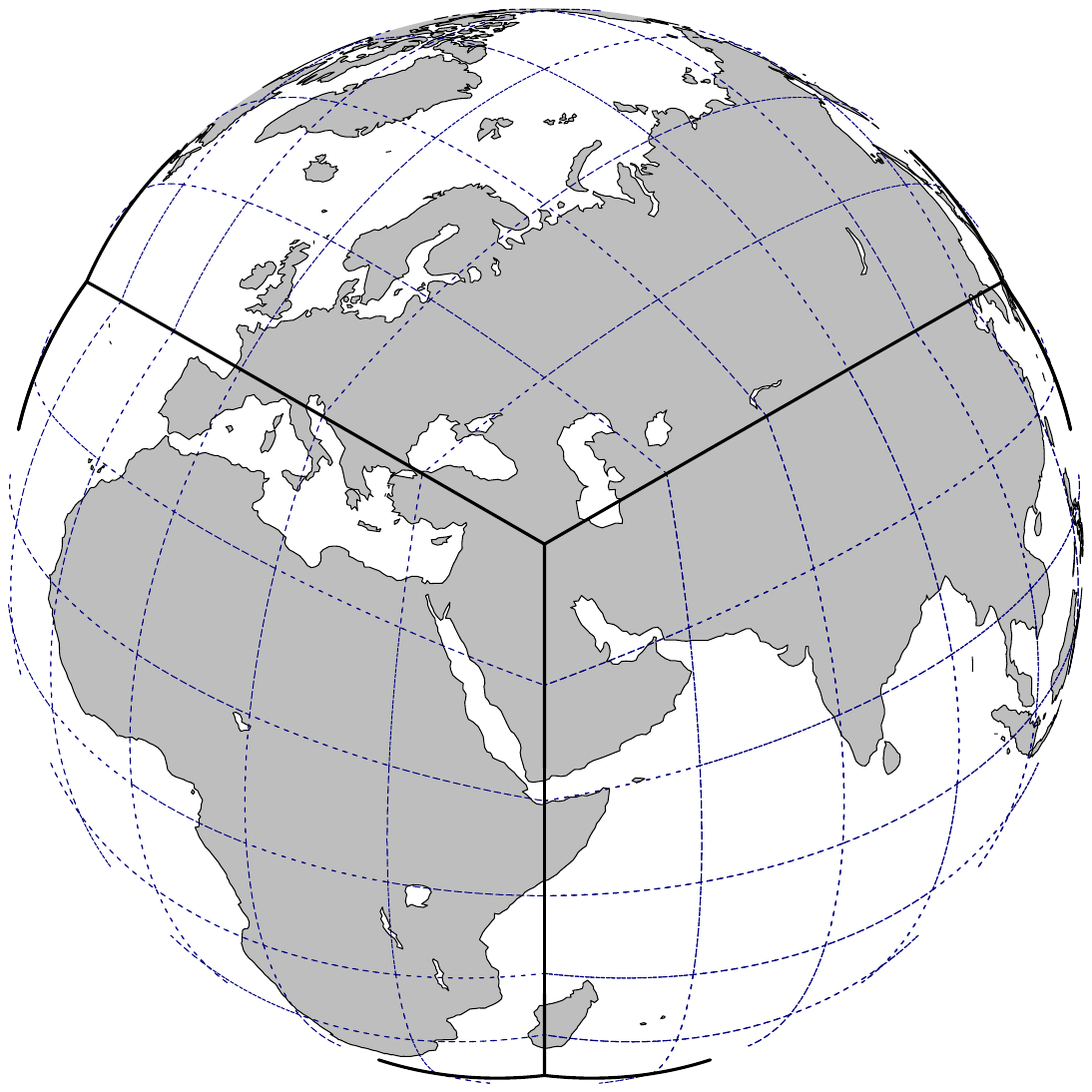}
    \includegraphics[trim={0 6cm 0 5cm},clip,width=0.99\linewidth]{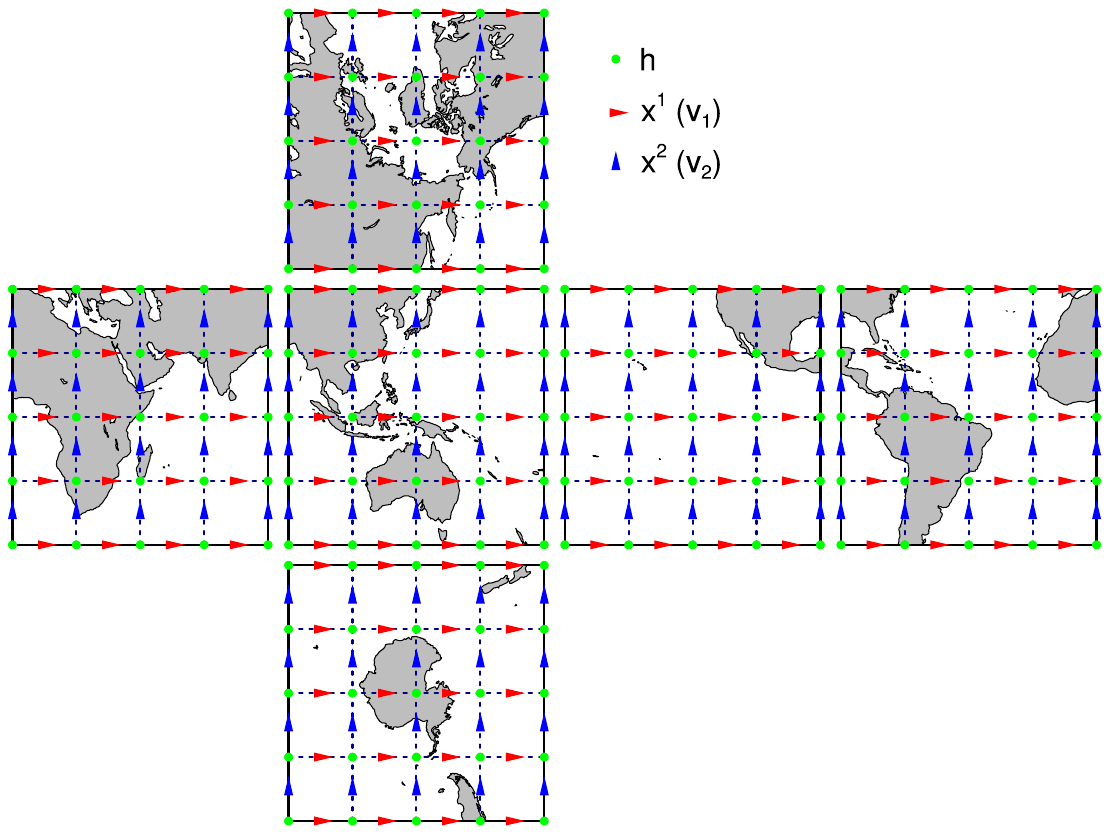}
    \caption{Arrangement of blocks and grid points in each block of equiangular gnomonic cubed-sphere grid.}
    \label{fig:grid}
\end{figure}

In this section, we show how SBP operations described in one spatial dimension may be extended to 2D multi-block logically rectangular grids via tensor (or Kronecker) products. 

First, consider closed multi-block computational domain with $N_b$ conforming logically-rectangular blocks with $(x_1, x_2)$ coordinates in each. Without loss of generality we will assume grid blocks to be squares of equal size. Each block is discretized using three types of points
\begin{alignat}{3}\label{eq:xh}
    \mathbf{x}^h_{ijp} = ({x}^v_{i}, {x}^v_{j}), \quad && i&=\overline{1,N+1}, \  && j = \overline{1,N+1}, \\ 
    \label{eq:x1}
    \mathbf{x}^1_{ijp} = ({x}^c_{i}, {x}^v_{j}), \quad && i&=\overline{1,N}, \ && j = \overline{1,N+1}, \\
    \label{eq:x2}
    \mathbf{x}^2_{ijp} = ({x}^v_{i}, {x}^c_{j}), \quad && i&=\overline{1,N+1}, \ && j = \overline{1,N},
\end{alignat} 
where $x_v$ and $x_c$ are defined according to (\ref{eq:xc_xv}) and $p\in[1, N_b]$ is a block index. The arrangement of blocks and points inside block is shown in Fig. \ref{fig:grid} using equiangular gnomonic cubed-sphere grid \cite{Rancic_CSSW_1996} as an example.

Column vectors
\begin{equation}
    \mathbf{h} \in \mathbb{R}^{N_b (N+1)^2\times1}, \quad \mathbf{v}_1, \mathbf{v}_2 \in \mathbb{R}^{N_b N (N+1)\times1},
    \end{equation}
contains grid functions values at the $\mathbf{x}^h$,  $\mathbf{x}^1$, $\mathbf{x}^2$ grids, respectively. The numbering of values in the column vectors above is first in the $x_1$ direction, then in $x_2$ direction and finally in the block index dimension, e.g. the number of $\mathbf h$ component corresponding to the value at $x^h_{ijp}$ point is $(p-1)(N+1)^2+(j-1)(N+1)+i$.

Next, we define the differentiation and interpolation SBP operators in 2D multi-block domain defined above using Kronecker products. The differentiation operators become
\begin{align}
    & D_{h1} = I_{N_b} \otimes I_{N+1} \otimes D_{vc}  \in \mathbb{R}^{N_1\times N_h}, \ 
    & D_{h2} = I_{N_b} \otimes D_{vc}  \otimes I_{N+1} \in \mathbb{R}^{N_2\times N_h}, \\
    & D_{1h} = I_{N_b} \otimes I_{N+1} \otimes D_{cv}  \in \mathbb{R}^{N_h\times N_1}, \ 
    & D_{2h} = I_{N_b} \otimes D_{cv} \otimes I_{N+1}  \in \mathbb{R}^{N_h\times N_2}, 
\end{align}
where $I_{N+1}$ and $I_{N_b}$ are the $(N+1)\times(N+1)$ and ${N_b}\times{N_b}$ identity matrices, respectively; $N_h = {N_b} (N+1)^2$, $N_1=N_2 = {N_b}(N+1)N$. Using these notations, $D_{h1}$ is the differentiation matrix that approximates $\partial/\partial x_1$ at $\mathbf{x}^1$ grid nodes using values at $\mathbf{x}^h$ grid nodes; $D_{h2}$ approximates $\partial/\partial x_2$ at $\mathbf{x}^2$ grid nodes using values at $\mathbf{x}^h$ grid; $D_{1h}$, $D_{2h}$ matrices approximate $\partial / \partial x_1$, $\partial / \partial x_2$ at $\mathbf{x}^h$ points using values at $\mathbf{x}^1$, $\mathbf{x}^2$ nodes, respectively.

Analogously, we introduce the quadrature matrices
\begin{equation}\label{eq:sbp_2d_norm}
     H_{h} = I_{N_b} \otimes H_v \otimes H_v \in \mathbb{R}^{N_h\times N_h}, \ 
     H_{1} = I_{N_b} \otimes H_v \otimes H_c \in \mathbb{R}^{N_1\times N_1}, \
     H_{2} = I_{N_b} \otimes H_c \otimes H_v \in \mathbb{R}^{N_2\times N_2},
\end{equation}
and interpolation operators
\begin{align}\label{eq:sbp_2d_interp1}
    & P_{h1} = I_{N_b} \otimes I_{N+1} \otimes P_{vc}  \in \mathbb{R}^{N_1\times N_h}, \ 
    & P_{h2} = I_{N_b} \otimes P_{vc}  \otimes I_{N+1} \in \mathbb{R}^{N_2\times N_h}, \\
    \label{eq:sbp_2d_interp2}
    & P_{1h} = I_{N_b} \otimes I_{N+1} \otimes P_{cv}  \in \mathbb{R}^{N_h\times N_1}, \ 
    & P_{2h} = I_{N_b} \otimes P_{cv} \otimes I_{N+1}  \in \mathbb{R}^{N_h\times N_2}.
\end{align}
Here $P_{h1}$, $P_{h2}$ interpolates values from $\mathbf{x}^h$ grid to $\mathbf{x}^1$ and $\mathbf{x}^2$ grids, respectively, and vice versa for the $P_{1h}$, $P_{2h}$ matrices.   

Using the properties of Kronecker product, we can derive the 2D analogue of the SBP property  (\ref{eq:sbp_property}):
\begin{equation}\label{eq:2d_diff_matrices}
\begin{gathered}
     H_1 D_{h1}  = - D_{1h}^T H_h + I_{N_b} \times I_{N+1} \times R_{cv} \equiv - D_{1h}^T H_h + R_{1h}^T \\
     H_2 D_{h2}  = - D_{2h}^T H_h + I_{N_b} \times R_{cv} \times I_{N+1} \equiv - D_{2h}^T H_h + R_{2h}^T,
\end{gathered}
\end{equation}
and 2D analogue of SBP preserving property (\ref{eq:sbp_interp_requirement}) for interpolation operators
\begin{equation}\label{eq:sbp_2d_interp_property}
    H_h P_{1h} = P_{h1}^T H_1, \ H_h P_{2h} = P_{h2}^T H_2.
\end{equation}

Further, we will use the following notations:

\begin{equation}\label{eq:Dvh_Dhv}
D_{\mathbf{v} h} =
\begin{bmatrix}
    D_{1h} & D_{2h}
\end{bmatrix}, \
D_{h\mathbf{v}} =
\begin{bmatrix}
    D_{h1} \\
    D_{h2}
\end{bmatrix},
\end{equation}

\begin{equation}
    \mathbf{v} = \begin{bmatrix}
        \mathbf{v}_1 \\
        \mathbf{v}_2
    \end{bmatrix}, \ 
    R = \begin{bmatrix}
        R_{1h}\\  R_{2h}
    \end{bmatrix}.
\end{equation}

From (\ref{eq:2d_diff_matrices}) it follows:
\begin{equation}
    \begin{gathered}
    \mathbf{v}^T H_v D_{h\mathbf{v}} \mathbf{h} + \mathbf{h}^T H_h D_{\mathbf{v} h} \mathbf{v} = \mathbf{h}^T R_{1h}\mathbf{v}_1 + \mathbf{h}^T R_{2h}\mathbf{v}_2 \equiv [\mathbf{h}^T \ \mathbf{h}^T] R\mathbf{v}.
    \end{gathered}
\end{equation}

The presented SBP operators provide the foundation for constructing governing equation approximations on curvilinear multi-block meshes. To incorporate interface conditions, similar to the one-dimensional case, we need to augment these operators with SAT corrections and projection operators. Expressions for these operators will be presented in the following sections.

\section{Test problem formulation and discretization} \label{sec:model formulation}
\subsection{Analytical setup}
The linearized shallow water equations (\ref{eq:shallow_water}) inside each block of the domain are written as
\begin{subequations}\label{eq:curvilinear_sw}
    \begin{equation}
        \frac{\partial v_1}{\partial t} = f J v^2-\frac{\partial h}{\partial x_1},
    \end{equation}
    \begin{equation}
        \frac{\partial v_2}{\partial t} =-f J v^1-\frac{\partial h}{\partial x_2},
    \end{equation}
    \begin{equation}
        \frac{\partial h}{\partial t} = -\frac{\mathcal{H}}{ J}\left(\frac{\partial J v^1}{\partial x_1}+
        \frac{\partial J v^2}{\partial x_2}\right),
    \end{equation}
\end{subequations}
where $v_1$ and $v_2$ are covariant flow velocity components in $x_1$ and $x_2$ directions respectively, $v^1$ and $v^2$ are contravariant components, $J$ is the curvilinear  mapping Jacobian. Eqs. (\ref{eq:curvilinear_sw}) accomplished with the continuity conditions at the domain blocks boundaries govern dynamics of the shallow water in the whole domain.

We use covariant velocity components $v_1$ and $v_2$ as prognostic variables and contravariant components are auxiliary: $(v^1, v^2)^T = \mathcal Q (v_1, v_2)^T$, where $\mathcal Q$ is the contravariant metric tensor.

\subsection{Spatial discretization}

For the spatial discretization of Eqs. (\ref{eq:curvilinear_sw}) we use grid defined in Sect. \ref{sec:2d_sbp_sat}. The placement of variables is the following: $h$ field is stored at $\mathbf x^h$ points, and the velocity components $v_1$ and $v_2$ are placed at $\mathbf x^1$ and $\mathbf x^2$ points, correspondingly (see Fig. \ref{fig:grid} for visual reference). This is basically Arakawa type C grid \cite{ArakawaLamb1977}, which is known to deliver optimal numerical dispersion properties for a variety of geophysical equation systems \cite{Randall_GeostAdj,ThuburnWoolings2005}. Then, using 2D SBP-operators defined in Sect. \ref{sec:2d_sbp_sat} we arrive to SBP-SAT-Projection discretization:
\begin{subequations}\label{eq:discrete_sw}
    \begin{equation}\label{eq:discrete_sw_momentum}
        \frac{\mathrm{d} \mathbf{v}}{\mathrm{d} t}
         =  F\tilde{\mathbf{v}}
            - g\, D_{h\mathbf{v}} A_h \mathbf h,
    \end{equation}
    \begin{equation}
        \frac{d \mathbf h}{d t} = -\mathcal{H} A_h J_h^{-1}(D_{\mathbf{v}h}+S) J_\mathbf{v} \tilde{\mathbf{v}},
    \end{equation}
    \begin{equation}
       \tilde{\mathbf{v}} = 
        Q \mathbf{v},
    \end{equation}
\end{subequations}
where $\mathbf{v} = (\mathbf v_1^T,\mathbf v_2^T)^T$ contains column-vectors representing $v_1$ and $v_2$ functions, $\tilde{\mathbf{v}} = ((\mathbf {v}^1)^T,(\mathbf v^2)^T)^T$ are the same for $v^1$ and $v^2$, $\mathbf h$ is the column-vector representing $h$ field, $J_h$ is the diagonal matrix with mapping Jacobian values at $\mathbf x^h$ grid points and $J_\mathbf v = \mathbf{diag}(J_1, J_2)$ is the block diagonal matrix with $J_1$, $J_2$ being diagonal matrices containing Jacobian values at $\mathbf x^1$, $\mathbf x^2$ points,  respectively, $F$ is the Coriolis terms operator,  $Q$ is the discrete contravariant metric tensor operator, $A_h$ and $S$ are the projection and SAT terms operators, respectively, which are used to impose grid-blocks interface conditions. The explicit form of operators $Q$, $F$, $S$, $A_h$ is given below.

\textbf{A comment on the placement of variables.} The placement of variables within staggered grid (Fig. \ref{fig:grid}) is motivated by our goal to impose grid-block interface conditions in a SAT-Projection manner presented in Sect. \ref{sect:sat_proj}. This approach reduces the stiffness and removes the spurious stationary modes associated with discontinuities at the grid-block interfaces (see Sect. \ref{sect:laplace_spectrum}). We selected $h$ to be the variable projected to the interface-continuous subspace, because projecting vector fields is more complicated. This ultimately leads to the placement of $h$-points to the grid cells vertices so that they are located at the boundaries between the grid blocks.

We impose no explicit continuity conditions on the velocity field. However, as we will show later (Sect. \ref{sect:mimetic:vt}), our choice of discrete operators result in continuity of the tangential velocity component at the grid blocks interfaces. Consequently, covariant velocity components are also continuous at the grid block boundaries (if the possible change of local basis orientations between blocks is taken into account).

The grid we utilize is $\Delta x/2$ shifted from the more conventional setup with $h$-points placed at the cell centers $\mathbf x_{ijp} = (x_i^c, x_j^c)$ (e.g., as used in \cite{Ullrich_SW_2010,ShashkinGoymanExpSL,Melvin_mixedfe_2019}). Using C-staggering with $h$-points at the cell centers leads to $v_1$ and $v_2$ components placement at $\mathbf x^2$ and $\mathbf x^1$ points, respectively. Then, trying to use SAT-Projection interface conditions with such a grid one arrives to the need to project massflux vector components $Jv^1$ and $Jv^2$ to the edge continuous subspace. For energy conservation this projection is required to be orthogonal and thus application of $H_\mathbf v J_\mathbf v Q$ inverse is needed (that we consider prohibitive for the practical numerical scheme).

\subsection{Projection and SAT-term operator}

The projection matrix $A_h \in \mathbb{R}^{N_h \times N_h}$ ensures grid functions interface continuity at $h$-points, while the SAT-term matrix $S \in \mathbb{R}^{N_h \times (N_1+N_2)}$ provides correction proportional to the flux mean at the interfaces, analogously to 1D case (\ref{eq:1d_sat_vc}). The particular form of these matrices is determined by the topology of the multi-block grid utilized (connections between the blocks, their mutual orientation), hence, we will define these operators in terms of their action on the grid functions at interface locations.

Let $m$ be the index of $\mathbf{h}$ element corresponding to some non-corner interface point on a grid block $p$, and $m^*$ be the index of $\mathbf{h}$-element corresponding to the physically same point but from the neighbour block, then
\begin{equation}\label{eq:Ah_proj}
    \left(A_h \mathbf{h}\right)_m = \frac{(J_hH_h \mathbf{h})_{m} + (J_hH_h \mathbf{h})_{m^*}}{(J_hH_h)_{mm}+ (J_hH_h)_{m^* m^*}}.
\end{equation}
Application of $A_h$ matrix results in $(A_h\mathbf{h})_m = (A_h\mathbf{h})_{m^*}$.
For corner points, it is necessary to consider the contribution of all neighboring grid blocks:
\begin{equation}\label{eq:Ah_proj_vertex}
    \left(A_h \mathbf{h}\right)_m = \frac{(J_h H_h \mathbf{h})_{m} + \sum_{m^*} (J_h H_h \mathbf{h})_{m^*}}{(J_h H_h)_{mm}+ \sum_{m^*} (J_h H_h)_{m^* m^*}}.
\end{equation}
Since $A_h$ is an orthogonal projector matrix, it follows that $A_h^2 = A_h$ and $(J_h H_h A_h)^T =J_h H_h A_h$.

The $S$ matrix can be written as follows:
\begin{equation}\label{eq:SAT-matrix}
    S = H_h^{-1}
    \begin{bmatrix}
        I_{N_h} & I_{N_h}
    \end{bmatrix}
    A_n R.
\end{equation} 
The $A_n$ matrix approximates mean value of the interface fluxes acting on the vector
\begin{equation}
    \mathbf{u} = J_{\mathbf{v}}\tilde{\mathbf{v}} = 
    \begin{bmatrix}
        \mathbf{u}^1 \\ \mathbf{u}^2
    \end{bmatrix}.
\end{equation}
Assume $\left(R\mathbf{u}\right)_m$ element contains extrapolated value of the $\mathbf{u}$ vector component at some interface point of a grid block, and $\left(R\mathbf{u}\right)_{m^*}$ element corresponds to the extrapolated value at the physically same point of a neighbor block interface (in the case of corner points, the vector component in parallel direction is used). The matrix $A_n$ modifies these values according to:
\begin{equation}
    (A_n R\mathbf{u})_m = (A_n R\mathbf{u})_{m^*} = \frac{1}{2}\left((R\mathbf{u})_m + (R\mathbf{u})_{m^*}\right).
\end{equation}
For instance, consider the north interface point $x^h_{i(N+1)p}$. The interface flux at this point is determined by the extrapolated value of $\mathbf{u}^2$ elements which is stored at the $m=N_h+(p-1)(N+1)^2 + N(N+1) + i$ element of $R\mathbf{u}$. Next, suppose that the west interface of a block $p^*$ contains the same physical point at $x^h_{(N+1)j^*p^*}$. The flux at this point for block $p^*$ is determined by the extrapolated value of the $\mathbf{u}^1$ component contained at the element of the vector $R\mathbf{u}$ with the index $m^* = (p^*-1)(N+1)^2+(j^*-1)(N+1)+N+1$.

\subsection{Contravariant metric tensor operator}
The operator $Q$ is used to transform covariant components $\mathbf v_1$, $\mathbf v_2$ to contravariant components $\mathbf v^1$, $\mathbf v^2$ to calculate divergence and Coriolis force. Besides, this operator defines the dot product of two vector fields. The approximation for $Q$ is similar to \cite{OReilly_staggered_sbp_2020}:
\begin{equation} \label{eq:co2contra}
    \begin{bmatrix}
        \mathbf v^1 \\ \mathbf v^2
    \end{bmatrix}=
    \begin{bmatrix}
        Q^{11}_{1} & J^{-1}_{1} P_{h1} J_h Q^{12}_h P_{2h} \\
        J^{-1}_{2} P_{h2} J_h Q^{12}_h P_{1h} & Q^{22}_{2}
    \end{bmatrix}
    \begin{bmatrix}
        \mathbf v_1 \\ \mathbf v_2
    \end{bmatrix},
\end{equation}
where $Q^{lm}_q$ are diagonal matrices of analytic contravariant metric tensor component $lm$ at $\mathbf x^q$  grid points.

The approximation (\ref{eq:co2contra}) follows the ideas of \cite{OReilly_staggered_sbp_2020}: the contravariant metric tensor off-diagonal terms action is calculated at $h$-points. Thereby, velocity components are interpolated to the $h$-points, then multiplied by $J_h Q^{12}_h$ and the result is interpolated back to the velocity points. 

\subsection{Coriolis terms operator}

The Coriolis effect plays the crucial role in such phenomena as geostrophic adjustment and Rossby waves that dominate ocean and atmospheric dynamics. Moreover, the subtle details of Coriolis terms discretization at the staggered grid significantly influence the numerical dispersion of Rossby modes as shown in \cite{thuburn2007rossby}. Below, we present several options for Coriolis terms approximation so that one can select the most suitable one for the specific problem. All considered Coriolis terms approximations are energy neutral.

The basic energy neutral approximation can be constructed following the ideas of regular latitude-longitude discretization \cite{thuburn2004CgridCons}
\begin{equation}\label{eq:coriolis_ini}
F = \underbrace{P_{h\mathbf v}}_3 \underbrace{C}_2 
\underbrace{P_{\mathbf v h} J_{\mathbf v}}_1,
\end{equation}
where
\begin{equation}
    C = \begin{bmatrix}
        & \textrm{diag}(f) \\ -\textrm{diag}(f) &
    \end{bmatrix}
\end{equation}
with $\textrm{diag}(f)$ being diagonal matrices of Coriolis parameter values at $h$-points.

The above formula (\ref{eq:coriolis_ini}) implies three-stage algorithm: (1) contravariant velocity components are multiplied by $J$ and interpolated to $h$ points (as recommended by \cite{thuburn2007rossby} for the better Rossby-waves dispersion characteristics), (2) Coriolis force components are calculated at $h$ points, (3) Coriolis force components are interpolated back to the velocity points.

The approximation (\ref{eq:coriolis_ini}) produces the tendencies that are discontinuous across the grid block edges because no interaction between the grid blocks are introduced. That leads to the discontinuity of the tangential velocity at the boundaries between blocks. It is currently unclear whether contributions to this discontinuity from linear dynamics (as well as contributions from nonlinear terms and forcing for real simulations) can accumulate and spoil long-term solution. Also, discontinuity of Coriolis tendencies decouples degrees of freedom associated with covariant velocity components at the edge points and thus might lead to the presence of nonphysical wave modes.

Energy-neutral Coriolis terms operator that preserves the continuity of tangential velocity at the grid block boundaries is
\begin{equation}\label{eq:cori_full}
    F = 
    \frac{1}{2}
    P_{h\mathbf v}
    \left(
        V J_{\mathbf h}C+J_{\mathbf h}C \tilde V
    \right)
    J_{\mathbf h}^{-1}
    P_{\mathbf v h}J_{\mathbf v},
\end{equation}
where $J_{\mathbf h} = \textrm{diag}(J_h,J_h)$. The operator $V$ modifies the covariant vector components defined at $h$-points to make the vector field continuous at the grid block boundaries:
\begin{equation}\label{eq:A_cov}
    V = Y \begin{bmatrix}A_h & & \\ & A_h & \\ & & A_h\end{bmatrix} X,
\end{equation}
where matrix $X$ transforms covariant vector components to Cartesian components and $Y$ transforms Cartesian components to the covariant ones (here we imply that the considered 2D domain is embedded into 3D space). The orthogonal projection (\ref{eq:Ah_proj}, \ref{eq:Ah_proj_vertex}) to scalar edge-continuous subspace is used for each of the three Cartesian components. The $\tilde V$ is the same as $V$ but for contravariant components:
\begin{equation}\label{eq:A_contra}
    \tilde V = X^T \begin{bmatrix}A_h & & \\ & A_h & \\ & & A_h\end{bmatrix} Y^T,
\end{equation}
where contravariant to Cartesian and Cartesian to contravariant components transforms are $Y^T$ and $X^T$, respectively, according to the basic curvilinear grid formalism.

The edge continuity of Coriolis operator (\ref{eq:cori_full}) requires that the same interpolation formulae are used for each pair of neighbour grid blocks along the shared boundary. This requirement is satisfied when the same SBP discretizations are used for each block. Then, application of $P_{h\mathbf v}$ does not break the edge continuity of Coriolis terms.

The Coriolis terms edge continuity at $h$ points itself is evident for the $V J_\mathbf h C$ part of expression \eqref{eq:cori_full}. For the $J_\mathbf h C \tilde V$ we can note that $J_\mathbf h C$ is just $f k\times$ applied to the vector at each point. This do not break the edge continuity provided both $k$ and $f$ are edge continuous (that is the case for the smooth domains such as sphere, ellipse, etc).

The approximation (\ref{eq:cori_full}) can be simplified if the special type of grid is used where the values of mapping Jacobian are continuous at grid-blocks boundaries, i.e. for each point shared by several blocks, the value of $J$ is the same for each sharing block. The example of such grid is the cubed-sphere grid \cite{Rancic_CSSW_1996}. In this case, $J_{\mathbf h}$ commutes with $V$ and $V C = C \tilde V$ and thus Eq. (\ref{eq:cori_full}) is equivalent to 
\begin{equation}\label{eq:coriolis_opt2}
F = P_{h\mathbf v}V C P_{\mathbf v h} J_{\mathbf v}.
\end{equation}
However, in the case of cubed-sphere we find that the following energy-neutral interface-continuous approximation is more accurate (see Sect. \ref{sect:num_exp} for more information on accuracy):
\begin{equation}\label{eq:coriolis_main}
F = J_{\mathbf v}^{-1} P_{h\mathbf v} V J^2_{\mathbf h}C P_{\mathbf v h}.
\end{equation}

\section{Mimetic properties}\label{sect:mimetic_prove}

In this section, we examine whether the presented scheme replicates the discrete analogues of key symmetry and conservation properties inherent to the analytical equations. These properties are often regarded as crucial in geophysical hydrodynamics.

\subsection{Correctness of the contravariant metric tensor}\label{sect:mimetic:Q}
Contravariant metric tensor determines discrete dot-product rule for vector fields:
\begin{equation}\label{eq:kinetic_energy_quad}
(\vec w,\vec v) =
\int_\Omega(w_1 v^1+w_2 v^2) J dx_1 dx_2 \approx 
\mathbf{w}^T H_\mathbf{v} J_\mathbf{v} \tilde{\mathbf{v}} \equiv 
\mathbf{w}^T H_\mathbf{v} J_\mathbf{v} Q {\mathbf{v}}
\end{equation}
where $H_\mathbf{v} = \mathbf{diag}(H_{1},H_{2})$. Symmetry and positive definiteness $(H_\mathbf{v} J_\mathbf{v} Q)^T = H_\mathbf{v} J_\mathbf v Q > 0$ are needed for this dot-product to be valid.

\textbf{Symmetry proof}. Using Eq. (\ref{eq:co2contra}):
\begin{equation}\label{eq:symmetryQ1}
J_\mathbf{v} H_\mathbf{v}  Q =   \begin{bmatrix}
        H_{1}J_{1} Q^{11}_{1} & H_1 P_{h1} Q^{12}_h J_{h} P_{2h} \\
        H_2 P_{h2} Q^{12}_h J_h  P_{1h} & H_{2}J_{2} Q^{22}_{2}
    \end{bmatrix}.
\end{equation}
All $H_q$, $J_q$ and $Q^{lm}_q$ are diagonal matrices and thus commute. Therefore, diagonal blocks of (\ref{eq:symmetryQ1}) are symmetric $(H_q J_q Q^{ll}_q)^T = H_q J_q Q^{ll}_q$. For off-diagonal blocks we also use the interpolation operators properties (\ref{eq:sbp_2d_interp_property}):
\begin{multline}
    (H_1 P_{h1} J_h Q^{12}_h P_{2h})^T =  (P_{2h})^T Q^{12}_h J_h (H_1 P_{h1})^T =\\= H_2 P_{h2} H_h^{-1} Q^{12}_h J_h H_h P_{1h} = H_2 P_{h2} Q^{12}_h J_h P_{1h}, 
\end{multline}
which means the upper-right block of (\ref{eq:symmetryQ1}) is the transpose of lower left block. Thus, $J_\mathbf{v} H_\mathbf{v}  Q$ is symmetric.

\textbf{Positive definiteness proof.} We denote the blocks of the matrix $H_{\textbf{v}}J_{\textbf{v}}Q$ as follows:
\begin{equation}
    H_{\textbf{v}}J_{\textbf{v}}Q = 
    \begin{bmatrix}
        H_{1}J_{1} Q^{11}_{1} & H_1 P_{h1} Q^{12}_h J_{h} P_{2h} \\
        H_2 P_{h2} Q^{12}_h J_h  P_{1h} & H_{2}J_{2} Q^{22}_{2}
    \end{bmatrix}
    \equiv
    \begin{bmatrix}
        W_{11} & W_{12}\\ W_{12}^T & W_{22}
    \end{bmatrix}
\end{equation}

The matrix's positive definiteness is equivalent to the positive definiteness of both $W_{11}$ and the Schur complement $W_{22} - W_{12}^TW_{11}^{-1}W_{12}$ \cite{horn2012matrix}. The elements of diagonal matrices
$W_{11}$ and $W_{22}$ are positive, and hence, $W_{11}, W_{22} >0$. Moreover, $W_{12}^T W_{11}^{-1} W_{12} \equiv (W_{11}^{-\frac{1}{2}} W_{12})^T (W_{11}^{-\frac{1}{2}}W_{12}) \ge 0$.
Therefore, the positive definiteness of $W_{22} - W_{12}^T W_{11}^{-1} W_{12}$ is equivalent \cite{horn2012matrix} to the condition
\begin{equation}\label{eq:pos_def_crit}
\rho(W_{22}^{-\frac{1}{2}}W_{12}^TW_{11}^{-1}W_{12}W_{22}^{-\frac{1}{2}})<1,
\end{equation}
where $\rho$ denotes the spectral radius.

To illustrate the essence of this condition, consider a grid with a linear mapping where the angle $\alpha$ between the directions of the covariant basis vectors is constant. The analytical expressions for the contravariant metric tensor and the mapping Jacobian are:
\begin{equation}
    \mathcal{Q} = \frac{1}{\sin^2\alpha}\begin{bmatrix}
        1 & -\cos\alpha\\
       -\cos\alpha & 1
    \end{bmatrix}, \ J = \sin\alpha.
\end{equation}
Since these characteristics are independent of grid points, substituting the expressions for the matrices $W_{lm}$ into the criterion (\ref{eq:pos_def_crit}) yields:
\begin{equation}\label{eq:criterion_const_skew}
\cos^2\alpha\cdot \rho\left(P_{1h}P_{h1}P_{2h}P_{h2}\right)=\left(\cos\alpha\cdot\rho\left(P_{cv}P_{vc}\right)\right)^2 < 1.
\end{equation}

This result indicates that the positive definiteness of the discrete contravariant metric tensor is determined by the grid's skewness. For the interpolation operators derived in this work, $\rho\left(P_{cv}P_{vc}\right) \le 1.22$, which implies that  $\alpha > \arccos{(1/1.22)} \approx 35^{\circ}$. This angle provides a rough estimate 
of the maximum skewness that a general curvilinear grid can have while ensuring the scheme remains stable. 

A rigorous test of positive definiteness involves calculating the maximum eigenvalue of the sparse symmetric positive semidefinite matrix. This can be performed independently for each grid block, thanks to the Kronecker-product structure of the corresponding matrix. 

The computation of the matrix's spectral radius on the cubed-sphere grid shows that the criterion \eqref{eq:pos_def_crit} is satisfied with a significant margin. For all interpolation operators on grids with $N\le 2048$, the spectral radius $\rho < 0.37$. This result agrees well with the constant skewness estimate \eqref{eq:criterion_const_skew}, where $(1.22\cdot\cos{(60^\circ)})^2=0.3721$ (the minimum angle on a cubed-sphere is $60^\circ$).

\subsection{Gradient is anti-adjoint of divergence}\label{sect:mimetic:divgrad}

In this section, we demonstrate that the approximations of the divergence and gradient operators, as defined in System (\ref{eq:discrete_sw}), yield a discrete analogue of the Ostrogradsky-Gauss theorem for closed domains:
\begin{equation}
    \int_\Omega \vec v \nabla h d\Omega = -\int_\Omega h \nabla\cdot \vec v d\Omega,
\end{equation}
which indicates that the discrete gradient is the anti-adjoint of the divergence operator.

Using the discrete gradient and divergence operators from System (\ref{eq:discrete_sw}) and incorporating the vector fields quadrature defined in Equation (\ref{eq:kinetic_energy_quad}), we can express a discrete analogue of the Ostrogradsky-Gauss theorem:

\begin{equation}\label{eq:divgrad}
     (Q\mathbf{v})^T H_\mathbf{v} J_\mathbf{v} D_{h\mathbf{v}} A_h \mathbf{h} +
    \mathbf h^T H_h J_h A_h J_h^{-1} (D_{\mathbf vh }+S)
    J_v Q \mathbf{v} = 0.
\end{equation}

To prove this identity, we first use the symmetry of the projection operator $(H_h J_h A_h)^T = H_h J_h A_h$:
\begin{equation}
    (A_h \mathbf h)^T\left[ H_h (D_{\mathbf{v}h}+S) +(H_\mathbf{v} D_{h\mathbf{v}})^T \right] J_\mathbf{v}Q\mathbf{v} \equiv \tilde{\mathbf{h}} \left[ H_h (D_{\mathbf{v}h}+S) +(H_\mathbf{v} D_{h\mathbf{v}})^T \right] \mathbf{u}.
\end{equation}
Next, using the 2D SBP property (\ref{eq:2d_diff_matrices}) and expression for the SAT-term matrix (\ref{eq:SAT-matrix}), we can write:
\begin{equation}
\begin{gathered}
\begin{bmatrix}
    \tilde{\mathbf{h}}^T &
    \tilde{\mathbf{h}}^T
\end{bmatrix}
\left(R-A_n R\right) \mathbf{u} = \\
\sum_m \tilde{h}_m \left((R\mathbf{u})_m-\frac{1}{2} \left((R\mathbf{u})_m + (R\mathbf{u})_{m^*}\right)\right) = \\
\frac{1}{2}\sum_m \tilde{h}_m( (R\mathbf{u})_m- (R\mathbf{u})_{m^*}) = \frac{1}{2}\sum_m (R\mathbf{u})_m (\tilde{h}_{m}- \tilde{h}_{m^*}) = 0,
\end{gathered}
\end{equation}
where in the last transition we use edge continuity of $\tilde{\mathbf{h}} = A_h{\mathbf{h}}$ field i.e. $\tilde h_m = \tilde h_{m*}$.

\subsection{Mass-conservation}\label{subsec:mass_cons}
The discrete analogue for total mass is $\mathcal M = \mathbf 1_h^T H_h J_h \mathbf h$, where $\mathbf 1_h \in \mathbb{R}^{N_h \times 1}$ is the column vector of ones. The time-derivative of $\mathcal M$ writes as
\begin{multline}
    \frac{d \mathcal M}{d t} = \mathbf 1_h^T H_h J_h\frac{\partial \mathbf h}{\partial t} = 
    -\mathbf 1_h^T H_h J_h A_h J_h^{-1} (D_{\mathbf vh}+S) J_{\mathbf v} Q\mathbf{v}
    \overset{(*)}{=} \\
    =(Q \mathbf v)^T H_{\mathbf v} J_{\mathbf v} 
    D_{h \mathbf{v}} A_h \mathbf 1_h \overset{(**)}{=} 0,
\end{multline}
where the $(*)$ equality is due to the anti-adjoint property of divergence and gradient operators (Sect. \ref{sect:mimetic:divgrad}), $(**)$ is a result of all gradients used are exact (i.e. pointwise zero) for constant functions.

\subsection{Coriolis terms operator is energy neutral}\label{sec:mimetic:cori}
The analytic Coriolis force is pointwise kinetic energy neutral. Unlike the collocated grid case \cite{ShashkinGoymanTolstykh_SBP_SWE_2022}, where we prove the pointwise discrete analogue, with staggered grid we can achieve this only in the integral sense. Kinetic energy (\ref{eq:kinetic_energy_quad}) tendency exclusively due to Coriolis force action (\ref{eq:cori_full}) is
\begin{multline}\label{eq:cori_proof1}
    \left(\frac{d K}{d t}\right)_{cori} = \mathbf v^T J_\mathbf{v} H_\mathbf{v}  Q F  Q J_\mathbf v \mathbf v = \\
    \frac{1}{2}\mathbf u^T H_\mathbf v P_{h\mathbf v}\left(V J_{\mathbf h} C + J_{\mathbf h}C \tilde V \right)J_{\mathbf h}^{-1} P_{\mathbf v h} \mathbf u = \\
    =\frac{1}{2}\mathbf u_h^T \begin{bmatrix}
        H_h & \\ & H_h
    \end{bmatrix}
    \left(V J_{\mathbf h} C + J_{\mathbf h}C \tilde V \right)J_{\mathbf h}^{-1} \mathbf u_h,
\end{multline}
where $\mathbf u$ is $J_\mathbf v Q \mathbf v$ and $\mathbf u_h$ is the same interpolated to $h$-points and we use SBP-interpolation property (\ref{eq:sbp_2d_interp_property}). Denoting $H_{\mathbf h} = \textrm{diag}(H_h,H_h)$ we can proceed:
\begin{equation}\label{eq:cori_proof_intermediate}
     \left(\frac{d K}{d t}\right)_{cori} =
     \frac{1}{2}
     \left(\mathbf u_h^T J_{\mathbf h}^{-1}\right)
     H_{\mathbf h} J_{\mathbf h}
    \left(V J_{\mathbf h} C + J_{\mathbf h}C \tilde V \right)\left(J_{\mathbf h}^{-1} \mathbf u_h\right).
\end{equation}

Then, from Eqs. (\ref{eq:A_cov}, \ref{eq:A_contra}) it follows $(H_{\mathbf h} J_{\mathbf h} V)^T = H_{\mathbf h} J_{\mathbf h} \tilde V$ (that we leave without detailed proof for the sake of brevity) and thus
\begin{equation}
    (H_{\mathbf h} J_{\mathbf h} V J_{\mathbf h} C)^T = -J_{\mathbf h}C J_{\mathbf h} H_{\mathbf h} \tilde V
\end{equation}
due to anti-symmetry of $C$. Finally, noting that $J_{\mathbf h}H_{\mathbf h}$ commutes with $C$ (because vector rotation commutes with multiplication by a scalar) we can conclude:
\begin{equation}
    \left(J_{\mathbf{h}} H_{\mathbf{h}} V J_{\mathbf{h}} C\right)^T = 
    -J_{\mathbf{h}} H_{\mathbf{h}} J_{\mathbf{h}}C \tilde V
\end{equation}
and from (\ref{eq:cori_proof_intermediate}) it follows $(d K/d t)_{cori} = 0$.

The proof for Coriolis terms operator (\ref{eq:coriolis_main}) is generally the same and one should use that $V$ commutes with $J_{\mathbf h}$ when mapping Jacobian is continuous across grid-block boundaries.

\subsection{Energy conservation}\label{sect:total_energy_cons}
The total energy for the system of equations (\ref{eq:shallow_water}) is
\begin{equation}
    \mathcal E = \int_\Omega \left(\mathcal{H} \frac{\vec v\cdot\vec v}{2}+
                                   \frac{g}{2}{h^2}\right) d\Omega,
\end{equation}
where the integration is carried over the whole closed domain $\Omega$. The discrete total energy is
\begin{equation}
    \mathcal E \approx E = \frac{\mathcal{H}}{2}\mathbf v^T H_{\mathbf v} J_\mathbf{v}  Q \mathbf v + \frac{g}{2} \mathbf h^T H_h J_h \mathbf h.
\end{equation}

The conservation (i.e. $\partial E/\partial t=0$) can be easily proved using gradient-divergence anti-adjoint property (\ref{eq:divgrad}) and the linear kinetic energy neutrality of Coriolis terms operator (Sect. \ref{sec:mimetic:cori}).

Due to positive definiteness of $H_\mathbf v J_\mathbf v Q$, the total energy is positive definite quadratic functional and thus its conservation guarantees stability.

\subsection{Grid block interfaces continuity of tangential velocity}\label{sect:mimetic:vt}
The velocity component tangential to the grid block edge is defined by the covariant velocity component at the points ($\mathbf x_1$ or $\mathbf x_2$) lying at this edge. Suppose that for each edge dividing two neighbouring grid blocks the mapping between grid and physical coordinates is the same at both blocks \footnote{This assumption should not be limiting for grids conforming at the block edges as introduced in Sect. \ref{sec:2d_sbp_sat}.} (with account for possible difference in local bases orientations). Then, the tangential velocity at some point at the considered edge is continuous if the covariant velocity component values from both blocks coincide (again, with the possible sign change).

Therefore, the tangential velocity component at the grid block edges is continuous (or uniquely defined) if two conditions are satisfied: (1) the initial conditions are continuous at the edge and (2) all $\mathbf v$ tendencies are continuous at the edge during all time of simulation. The first condition is automatically satisfied for the physical problems dealing with continuous fields. The second condition
depends on the formulation of RHS of Eq. (\ref{eq:discrete_sw_momentum}). The pressure (height) gradient term of Eq. (\ref{eq:discrete_sw_momentum}) is edge-continuous as long as the $h$ field is edge continuous (that is assured by the presented SBP SAT-projection discretization) and for each two neighbour grid blocks the same differentiation formulae are used along the shared boundary. This is naturally satisfied if the same SBP approximation is used for each grid block. Analogously, the Coriolis terms approximations (\ref{eq:cori_full}, \ref{eq:coriolis_opt2}, \ref{eq:coriolis_main}) are interface-continuous if the same interpolation formulae are used at both sides of each grid block boundary.

To summarize, the continuity of velocity component tangential to the grid block edges is guaranteed if the SBP-SAT-Projection discretization based on the same 1D operators (Sect. \ref{sec:sbp_sat_1d}) is used in each grid block. 

\subsection{Zero curl of gradient}
The discrete analogue of $k\cdot \nabla\times\nabla \equiv 0$ is considered of primary importance for geophysical hydrodynamics for the reasons discussed in \cite{Stan_grid_rev}. On the staggered grid, curl $k\cdot\nabla\times\vec v = J^{-1}(-\partial v_1/\partial x_2+\partial v_2/\partial x_1)$ is naturally defined at the grid points shifted from all $h$, $v_1$ and $v_2$: $\mathbf x^\zeta_{ijp} = (x^c_i, x^c_j)$ (compare to Eqs. (\ref{eq:xh})-(\ref{eq:x2})). The discrete curl operator can be defined as:
\begin{equation}
    Z = J_\zeta^{-1}[-D_{1\zeta} \ D_{2\zeta}],
\end{equation}
where $J_\zeta$ is the diagonal matrix of the mapping Jacobian at $\mathbf x^\zeta$ points, $D_{1\zeta} = I_{N_b} \otimes D_{vc} \otimes I_{N}$, $D_{2\zeta} = I_{N_b}\otimes I_{N} \otimes D_{vc}$ are derivative operators in $x_2$ and $x_1$ directions respectively.

The gradient operator in the momentum Eq. (\ref{eq:discrete_sw_momentum}) is $D_{h\mathbf v}A_h$ and $D_{h\mathbf v}$ is itself defined in Eq. (\ref{eq:Dvh_Dhv}). Therefore, we get for the curl operator applied to the gradient operator:
\begin{equation}\label{eq:curlgrad}
    Z D_{h\mathbf v} = J_\zeta^{-1}(-D_{1\zeta}D_{h1}+D_{2\zeta} D_{h2})A_h = 
    J_\zeta^{-1} I_{N_b}\otimes (-D_{vc}\otimes D_{vc} + D_{vc}\otimes D_{vc})A_h \equiv 0.
\end{equation}

\section{Numerical experiments}\label{sect:num_exp}

\subsection{Preliminary notes}

To examine the properties of the presented SBP-discretizations, we conduct several tests in spherical geometry using equiangular gnomonic cubed-sphere grid \cite{Rancic_CSSW_1996} as in Fig. \ref{fig:grid}. We use four such grids with dimensions of $N_c = 48, 64, 96, 192$ that results in $N_c+1$ $h$-points along each cube's edge. We test three options of our schemes utilizing 2/1, 4/2 and 6/3-order accurate 1D staggered SBP-derivative operators, denoted Ch21, Ch42 and Ch63, correspondingly. The free parameters of the 6/3-order accurate 1D SBP derivative operators are defined by minimization of the objective function (\ref{eq:wave_opt_63}). The main reason for selecting objective function (\ref{eq:wave_opt_63}) rather than (\ref{eq:poly_opt_63}) is much smaller spurious reflection of short waves at the grid blocks boundaries noted in the experiment from Sect. \ref{sect:num_exp:poorly_resolved}.

The results of numerical experiments presented below are obtained using edge-continuous Coriolis force terms approximation (\ref{eq:coriolis_main}). Using simpler option (\ref{eq:coriolis_opt2}) results in almost 2 times higher errors in the tests with Coriolis force for all schemes (not shown). We also test edge-discontinuous variant of (\ref{eq:coriolis_main}) with $V$ dropped. Surprisingly, this lead to up to 2 times smaller error norms than for edge-continuous case, despite the jump of the tangential velocity at the grid block edges is considerable (not shown and not analyzed).

The results of staggered Ch schemes are compared to the Ah63 scheme that uses 6/3-order accurate collocated grid 1D SBP derivative. Ah63 is linear variant of \cite{ShashkinGoymanTolstykh_SBP_SWE_2022} scheme. All Ch schemes and Ah63 scheme share the same $h$-points at the cubed-sphere that facilitates direct comparison.

The classical explicit 4-th order Runge-Kutta scheme is used for integration in time, and time-steps values used at each grid are presented in Table \ref{tabl:num_exp_dt}. The time-step sizes are chosen in a way to keep the gravity waves Courant–Friedrichs–Lewy number (CFL) the same for each grid within the specific experiment. The resulting CFL values are rather small, aiming to minimize the influence of time-stepping errors on the overall accuracy of solution.

\begin{table}[H]
    \centering
    \begin{tabular}{ccccc}
    \hline
    $N_c$          & 48 & 64 & 96 & 192 \\ \hline
    $\Delta t$,\,s & 600 & 450 & 300 & 150 \\ \hline
    \end{tabular}
    \caption{Time-step values used in numerical experiments on the cubed-sphere grids of $6\times N_c\times N_c$ dimensions.}\label{tabl:num_exp_dt}
\end{table}

\subsection{Gaussian hill tests}
In this experiment, we consider the propagation of (inertia-) gravity waves caused by the initial Gaussian-shaped perturbation of the fluid height field:
\begin{equation}
    h(t=0) = \exp(-16 r^2 /a^2),
\end{equation}
where $r$ is the great-circle arc distance from the center of the perturbation, $a$ is the sphere radius taken to coincide with Earth radius. The initial velocity field is zero. Mean fluid height $\mathcal{H} = g^{-1}\left(2\pi a / 5\text{days}\right)^2 \approx 900\,\text{m}$ giving gravity waves phase speed equal to $c=\sqrt{g\mathcal{H}} = 2\pi a / 5\text{days}\approx 90\,\text{m}/\text{s}$.

Three variants of the test are examined:
\begin{enumerate}
    \item The center of perturbation is located at geographical latitude and longitude of $\varphi=0$, $\lambda=\pi$, which corresponds to the center of one of cubed-sphere blocks, the Coriolis parameter is zero.

    \item{The same as first variant, but the perturbation center is at $\varphi=\arcsin\left(\sqrt{1/3}\right)$, $\lambda=\pi/4$} that corresponds to one of the cubed-sphere vertices.

    \item The same as the second variant, but the Coriolis force is not zero with constant $f=10^{-4}\,s^{-1}$.
\end{enumerate}

In the first two variants of the experiment, the initial perturbation causes the propagation of the circular waves. After $2.5$ days of the model time the waves reach the point at the sphere which is opposite to the initial perturbation center and interfere to form the localized feature resembling initial perturbation shape. Then, in the course of next 2.5 days the waves propagate back to the initial perturbation location and so forth.

The main feature of the third variant is the Rossby-adjusted stationary mode with the $h$-field amplitude of about $0.5$ located at the initial perturbation position. The inertia-gravity wave part of the solution is similar to the first two variants except for the smaller amplitude and greater dispersion due to Coriolis effect.

\begin{figure}
    \centering
    \includegraphics[trim={0 9cm 0 9cm},clip,width=0.95\linewidth]{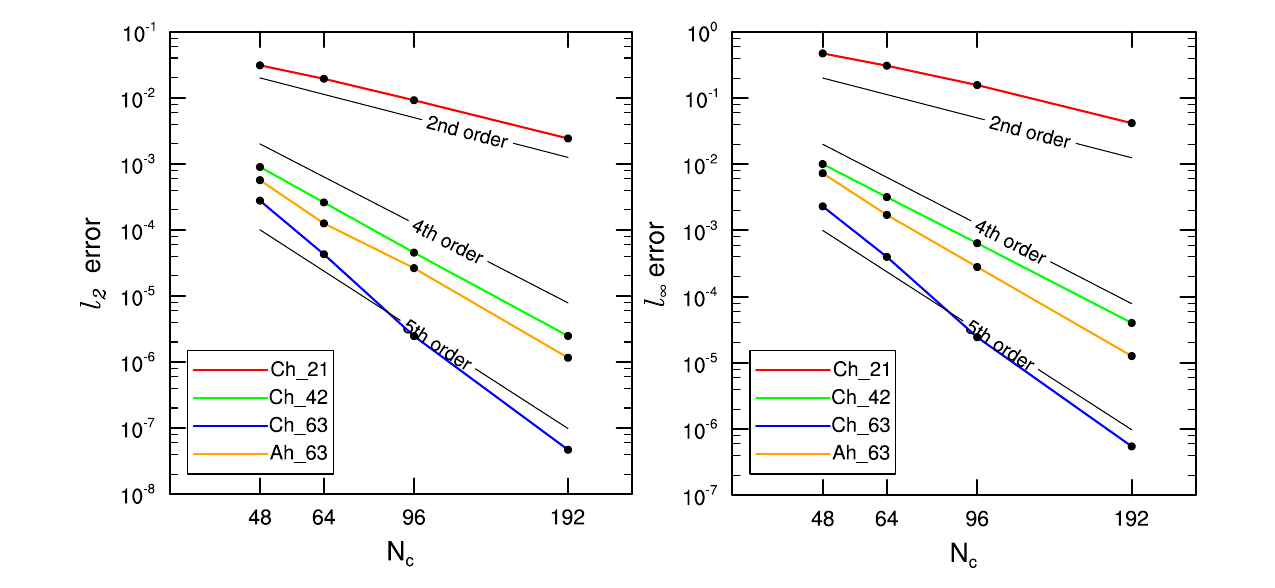}
    \caption{Dependence of 25-days simulation maximum values of $l_2$ and $l_\infty$ error norms for numerical solution of Gaussian hill test variant 1.}\label{fig:gauss_conv}
\end{figure}

\begin{figure}
    \centering
    \includegraphics[trim={0 9cm 0 9cm},clip,width=0.95\linewidth]{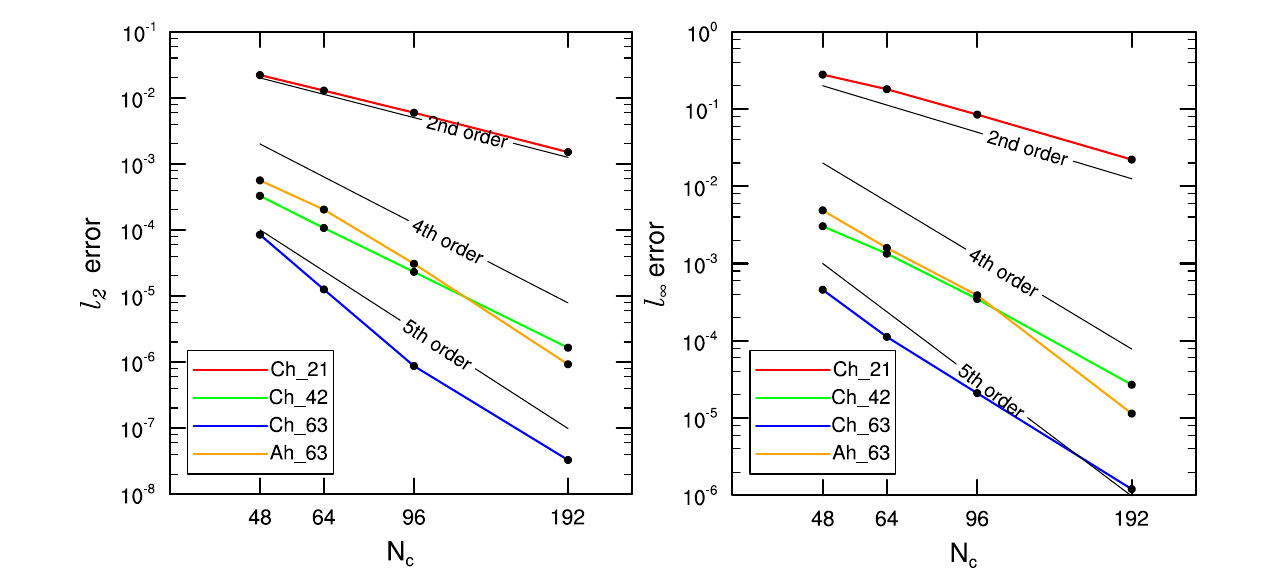}
    \caption{As Fig. \ref{fig:gauss_conv}, but for test variant 2.}\label{fig:gauss_corner_conv}
\end{figure}

\begin{figure}
    \centering
    \includegraphics[trim={0 9cm 0 9cm},clip,width=0.95\linewidth]{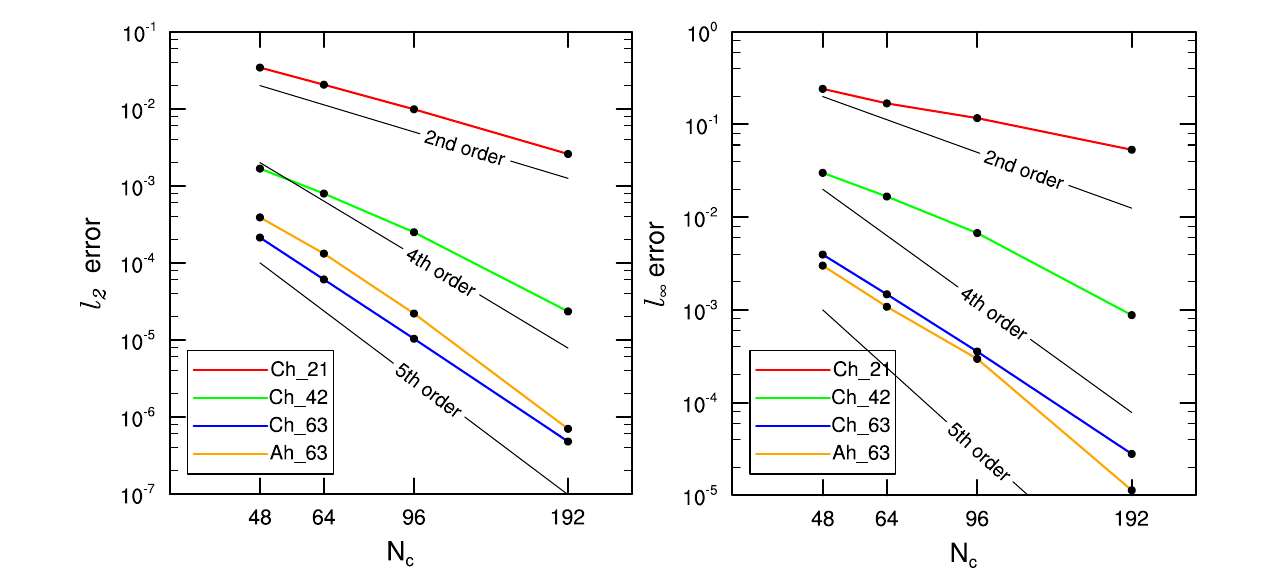}
    \caption{As Fig. \ref{fig:gauss_conv}, but for test variant 3.}\label{fig:gauss_fcorner_conv}
\end{figure}

The accuracy of the numerical solutions is measured using $l_2$ and $l_\infty$ error norms of $h$-field deviation from spectral transform model solution (which is exact up to double precision floating point arithmetics errors). The resolution dependence of the error norms maximums during 25 days of simulation is presented in Fig. \ref{fig:gauss_conv}, \ref{fig:gauss_corner_conv} and \ref{fig:gauss_fcorner_conv}. The corresponding solution convergence rates are listed in Table \ref{tabl:gauss_conv}.

\begin{table}[H]
    \centering
    \begin{tabular}{|c|c|c|c|c|}
    \hline
      & Ch21 & Ch42 & Ch63 &  Ah63 \\ \hline
    Gaussian hill 1 & 1.85 / 1.76 & 4.25 / 3.98 & 6.28 / 6.05 & 4.40 / 4.56 \\ \hline
    Gaussian hill 2 & 1.93 / 1.84 & 3.81 / 3.43 & 5.65 / 4.25 & 4.69 / 4.36 \\ \hline
    Gaussian hill 3 & 1.86 / 1.08 & 3.10 / 2.56 & 4.40 / 3.58 & 4.60 / 4.01 \\ \hline
    Solid rotation  & 2.14 / 1.43 & 3.33 / 2.70 & 3.72 / 3.41 & 4.13 / 3.95 \\ \hline 
    \end{tabular}
    \caption{Convergence rates in $l_2$ / $l_\infty$ norms for numerical solutions in the Gaussian hill  and solid rotation tests}\label{tabl:gauss_conv}
\end{table}

\begin{figure}
    \centering
    \includegraphics[trim={1cm 7cm 1cm 7cm},clip,width=\linewidth]{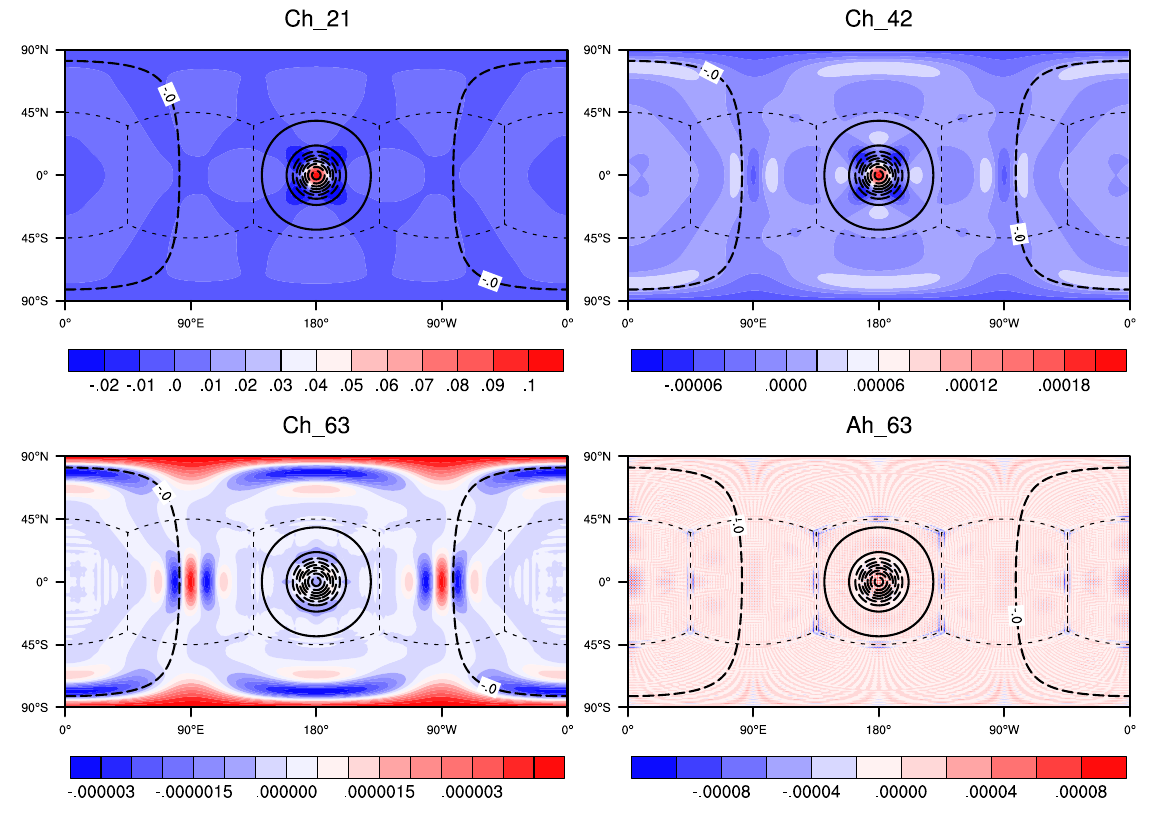}
    \caption{Solution (contours) and error (shading) for $h$-field at day 25 of Gaussian hill experiment variant 1, $N_c=96$. Dashed lines show blocks of the cubed-sphere grid.}
    \label{fig:gauss_map}
\end{figure}

\begin{figure}
    \centering
    \includegraphics[trim={1cm 7cm 1cm 7cm},clip,width=\linewidth]{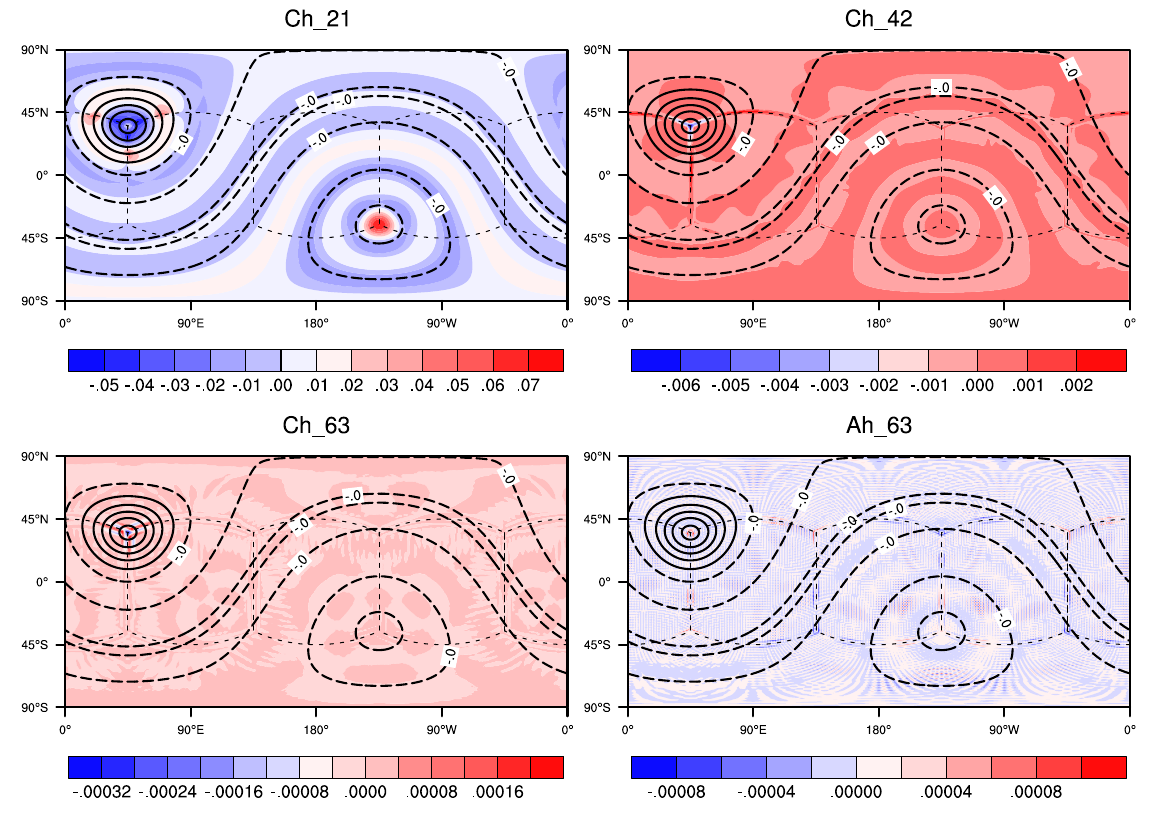}
    \caption{The same as Fig. \ref{fig:gauss_map}, but for Gaussian hill experiment variant 3.}
    \label{fig:gauss_map_fcorner}
\end{figure}

\textbf{Variant 1.} Staggered schemes solution convergence in both $l_2$ and $l_\infty$ norms is close to the scheme approximation order at grid block internal points (2nd for Ch21, 4th for Ch42 and 6th for Ch63). It seems that approximation error near the grid blocks boundaries does not significantly influence the solution, because the waves quickly pass this zone. However, this is not the case for the Ah63 scheme, which has an actual convergence order of about $4.5$ in both $l_2$ and $l_\infty$. The error pattern for Ch schemes is dominated by waves with spatial scales similar to those of the initial perturbation, as shown in Fig. \ref{fig:gauss_map}. Contrary, for Ah63 scheme, error fields are dominated by near grid-scale waves filling the domain as well as grid-scale features trapped near cubed-sphere edges.

The error norms for Ch63 solution are smaller than for Ah63 solution, and Ch63 convergence rate is higher. Therefore, using staggered grid benefits to the accuracy of reproduction of pure wave dynamics as compared to the case of collocated grid.

\textbf{Variant 2.} Convergence rates for Ch42 and Ch63 schemes drop down as compared to the test variant 1 (Table \ref{tabl:gauss_conv}). We believe SBP-interpolations needed for the calculation of contravariant velocity components and then divergence to be the primary source of error here. The collocated Ah63 scheme that do not need interpolation of velocity to calculate contravariant components shows convergence rates as in variant 1.

As in experiment variant 1, staggered grid Ch63 scheme outperforms collocated Ah63 in terms of $l_2$ and $l_\infty$ error norms and $l_2$ convergence rate. Ch63 and Ah63 $l_\infty$ convergence rates are quite similar. For the coarse grids lower-order Ch42 is more accurate than Ah63.

\textbf{Variant 3.} This is the most difficult case for the considered schemes because of the geostrophicaly-balanced feature persisting near one of the cube's vertices. To accurately maintain this feature, a numerical scheme needs to reproduce non-trivial balance of pressure (height) gradient and Coriolis forces and zero divergence of pure rotational velocity field. This is particularly challenging in the region of the grid with the most significant non-orthogonality and largest approximation errors of SBP difference and interpolation operators. Therefore, it is not surprising that the solution convergence rates for $2s$/$s$-order accurate Ch scheme drops to $s$ for $l_\infty$ and $s+1$ for $l_2$ (Ch63 is however slightly better with $s+1.4$ and $s+0.5$ convergence rates, respectively). The $l_2$ norm convergence rate of the Ah63 scheme does not drop down as compared to experiment variants 1 and 2, $l_\infty$ norm convergence rate drops down only slightly. We attribute this to the trivial contravariant metric tensor and Coriolis terms operators formulations at collocated grid that allows to circumvent the primary sources of error in this test.

In this test variant Ch63 is still more accurate than Ah63 in terms of $l_2$ norm and have similar values for $l_\infty$ error norm at all grids except the finest considered one. The $h$-field error for Ch schemes is dominated by the error near the cubed-sphere vertex. Ah63 error field is again combination of grid-scale waves and cubed-sphere edges trapped modes.

\subsection{Solid rotation test}

This is a linear variant of the standard meteorological shallow water test (number 2 from \cite{Williamson_tests}). The velocity field is the rotation around Earth geographical axis, $h$ field is in geostrophic balance with the wind field. With the non-constant Coriolis parameter $f = 2\Omega\sin\varphi$ this results in:
\begin{equation}
u = u_0 \cos\varphi',
\end{equation}
\begin{equation}
    h = h_0-\frac{a\Omega}{g}u_0 \sin^2\varphi',
\end{equation}
where $\varphi'$ is the latitude in the rotated geographical coordinate system such that the true North pole is located at $(\varphi, \lambda) = (\pi/4, 0)$.

In difference with the original (non-linear) test the height-field term proportional to $u_0^2$ is dropped because of metric terms in the velocity equation (\ref{eq:shallow_water}) are absent due to their non-linearity and do not need to be balanced by pressure gradient. The background value of fluid height $h_0= 29400/g$ from the original test setup is also taken as $\mathcal{H}$ in our model.

\begin{figure}
    \centering
    \includegraphics[trim={0 9cm 0 9cm},clip,width=0.95\linewidth]{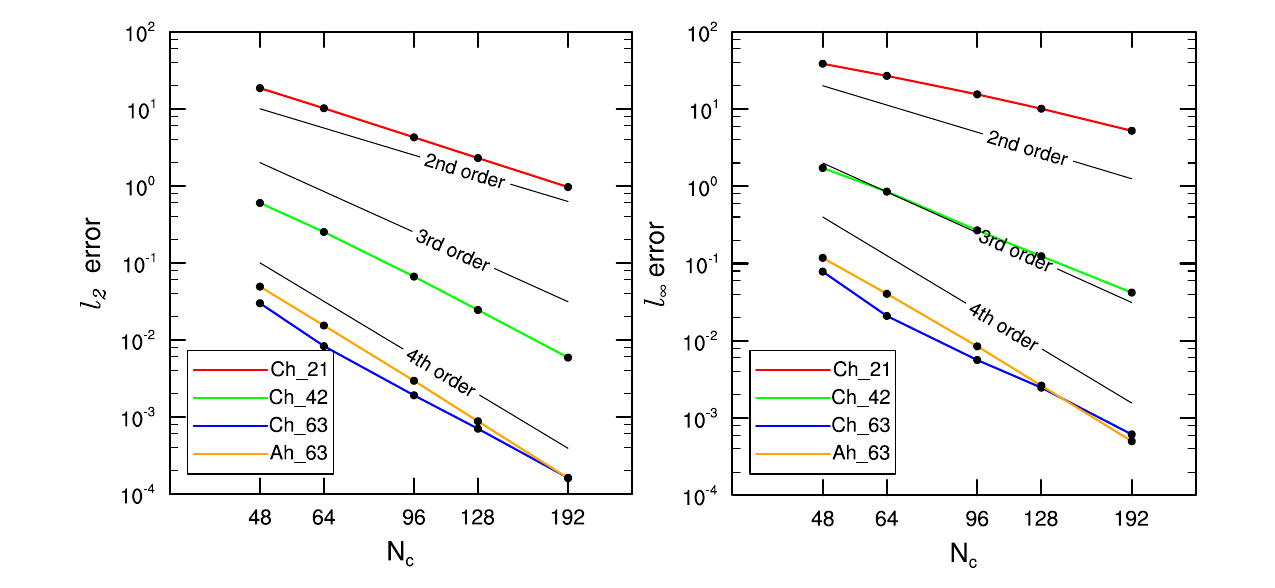}
    \caption{The values of $l_2$ and $l_\infty$ error norms for numerical solution of solid rotation test at day 10.}\label{fig:wt2_conv}
\end{figure}

\begin{figure}
    \centering
    \includegraphics[trim={1cm 7cm 1cm 7cm},clip,width=\linewidth]{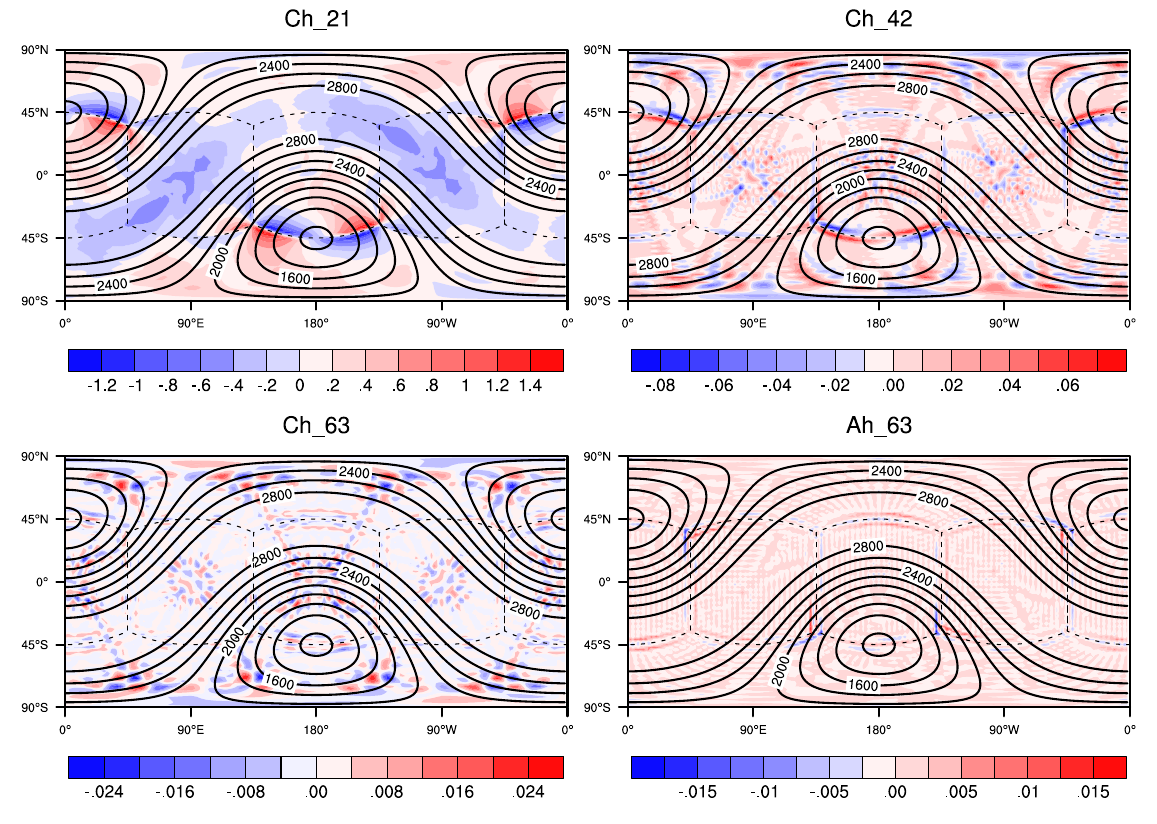}
    \caption{Solution (contours) and error (shading) for $h$-field at day 10 of solid rotation test, $N_c=48$. Dashed lines show blocks of the cubed-sphere grid.}
    \label{fig:wt2_map}
\end{figure}

The exact solution to this case is stationary, thus the numerical solution can be compared to initial conditions. The $h$-field $l_2$ and $l_\infty$ error norm values after 10 days of model time are shown in Fig. \ref{fig:wt2_conv}, the corresponding convergence rates are presented in Table \ref{tabl:gauss_conv}. The $l_2$ norm convergence rate for $2s$/$s$ order scheme (including Ah63) is near $s+1$ that is consistent with theoretical results from \cite{gustafsson1975convergence}. The $l_\infty$ convergence rate is significantly greater than $s$ for all Ch schemes. Fig. \ref{fig:wt2_map} confirms the difference between staggered and collocated schemes noted in Gaussian hill tests, i.e. the error pattern of Ah63 scheme is near grid-scale whereas Ch schemes error field is rather smooth and large scale.

\subsection{Poorly resolved waves dynamics}\label{sect:num_exp:poorly_resolved}
This test is similar to Gaussian hill test variant 1 but with Coriolis force turned on ($f=10^{-4}\,\text{s}^{-1}$) and that is more important, the initial perturbation is modulated by high frequency wave:
\begin{equation}
    h(t=0) = \exp(-16 r^2 /a^2)\cos^2(\nu\lambda)\cos^2(\nu\varphi),
\end{equation}
where $\nu=32,64$ that corresponds to the approximately $4\Delta x$ and $2\Delta x$ wavelength modulation at the $N_c=64$ cubed-sphere grid. In this experiment we test if the numerical schemes can deal with the features mimicking the external forcing (e.g. convective processes) to the atmospheric model which is usually high-frequency field with smooth envelope. 

Similarly to the Gaussian hill variant 3, the solution to this case is combination of inertia-gravity waves and stationary feature. According to the Rossby adjustment problem theory, the smaller spatial scale part of the initial conditions is primarily transformed to the inertia-gravity waves propagating away from the initial perturbation. The larger-scale part of initial conditions mainly remain as the stationary feature at the initial location, but some portion is also transformed to the inertia-gravity waves. The characteristic scale splitting the spatial frequency domain into short and long waves is the Rossby deformation radius $\lambda_R = \sqrt{g\mathcal{H}} / f$ with the value of about $900\,\text{km}$ for the selected test parameters.

\begin{figure}
    \centering
    \includegraphics[trim={1cm 5cm 1cm 5cm},clip,width=\linewidth]{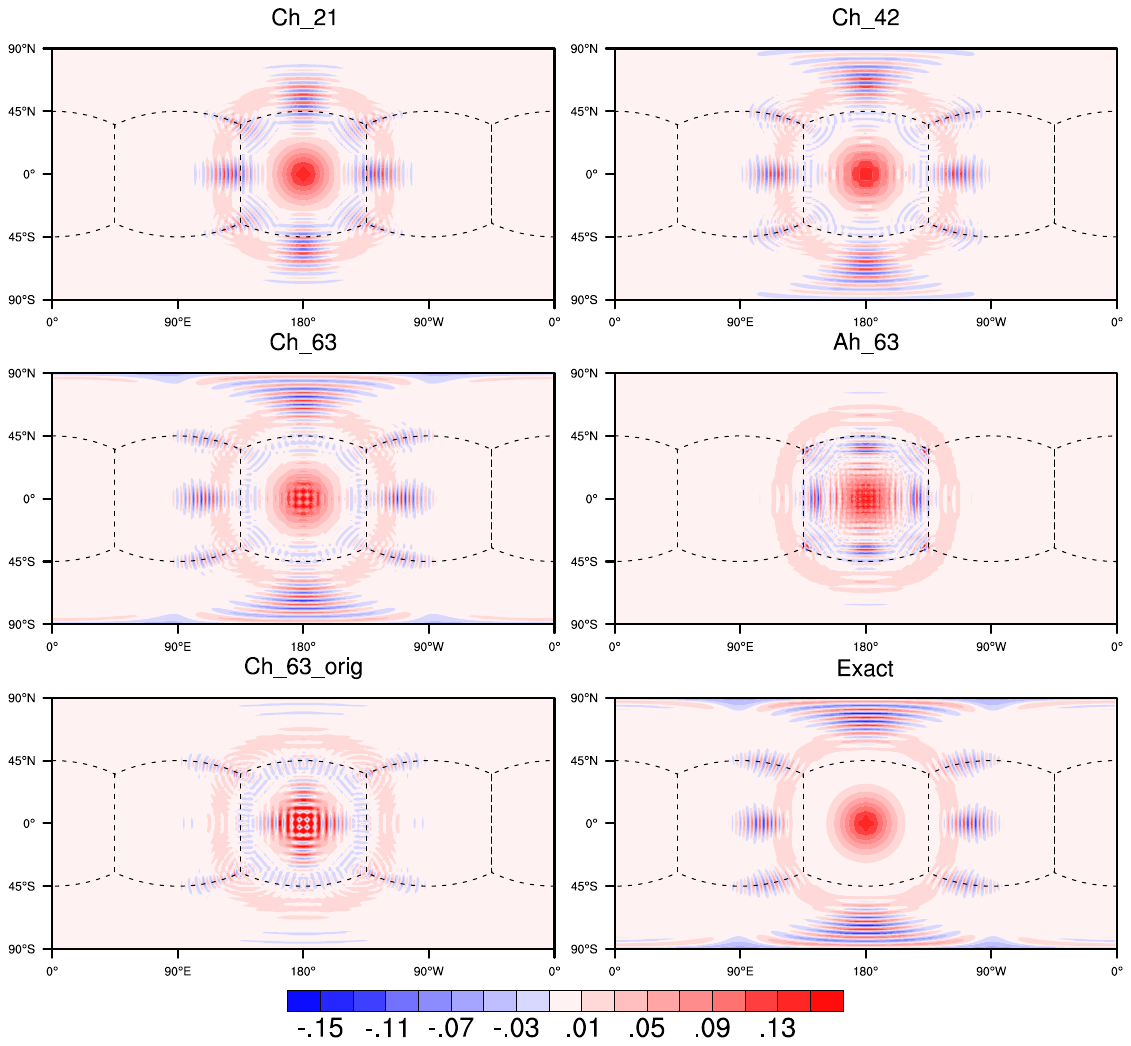}
    \caption{Solution for $h$-field after 25\,hours to poorly resolved waves dynamics test with $4\Delta x$ modulation.}
    \label{fig:bad_waves_4h}
\end{figure}

The numerical solutions for $4\Delta x$ modulation case at $t=25\,\text{hours}$ are compared to the exact one in Figure \ref{fig:bad_waves_4h}. We first analyze the effective propagation speed of $4\Delta x$ wave packets that can be estimated by the distance that the waves travel from the center of perturbation ($\varphi=0$, $\lambda=\pi$). One can see that the wave propagation speed increases as we increase the order of Ch schemes and for the Ch63 scheme it is close to the wave propagation speed in the exact solution (the wave packets crossing South/North poles are especially indicative). Second, one can note that all Ch scheme solutions suffer from the spurious waves reflection from the grid block interfaces. In the case of optimizing Ch63 scheme using polynomial objective function (\ref{eq:poly_opt_63}) ("Ch\_63\_orig" in Fig. \ref{fig:bad_waves_4h}) the $4\Delta x$ waves are almost fully reflected and do not propagate outside the grid block where the initial perturbation is located. Wave-accuracy optimization of 6/3-order 1D SBP derivative approximation (objective function \eqref{eq:wave_opt_63}) largely solves the reflection problem. The collocated grid Ah63 scheme solution is not only damaged by the wave reflection, but also suffers from the retarded propagation of short waves (slower than for Ch21) that is determined by the dispersion characteristics of collocated central differences.

The stationary mode of solution is often considered of primary importance in geophysical hydrodynamics, because (in more complicated model setup) it gives raise to slow Rossby-waves dynamics that dominates the atmospheric and oceanic flows. Therefore, it is particularly important for us to evaluate stationary mode that we approximately define as solution averaged over 600\,hours of simulation. With $4\Delta x$ modulation, all considered schemes including the collocated grid Ah63 scheme reproduce stationary mode close to the exact solution (not shown) despite the mentioned flaws in inertia-gravity wave propagation reproduction.

\begin{figure}
    \centering
    \includegraphics[trim={1cm 4cm 1cm 4cm},clip,width=\linewidth]{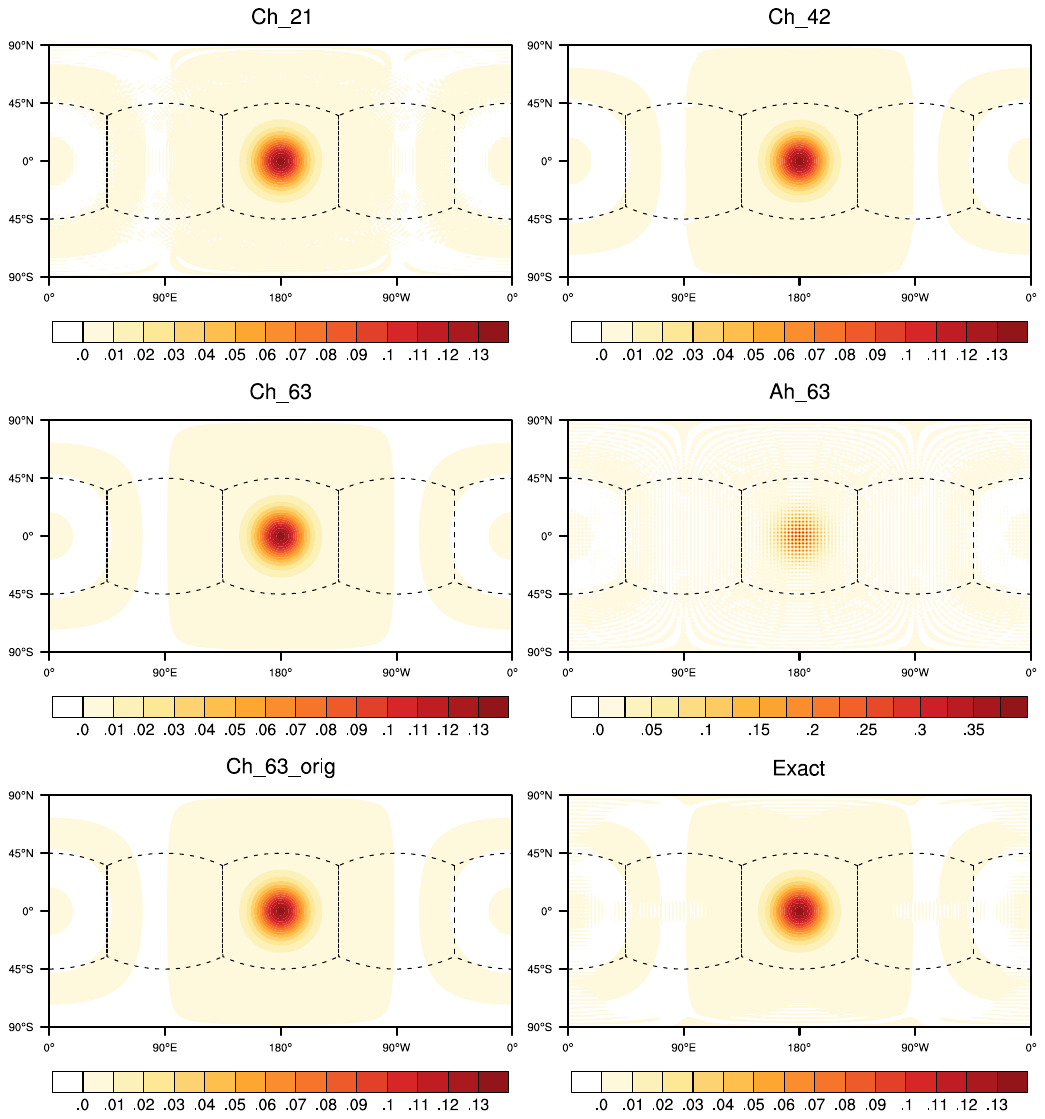}
    \caption{Stationary part (Averaging over $600$\,hours) of $h$-field solution to poorly resolved waves dynamics test with $2\Delta x$ modulation.}
    \label{fig:bad_waves_2h}
\end{figure}

In the $2\Delta x$ modulated case, the inertia gravity waves dynamics is not adequately resolved by the presented schemes. We observe retarded wave-packets propagation and spurious reflection from the grid blocks boundaries (not shown). Despite this, the stationary mode in the solution of staggered grid schemes is close to that of exact solution (see Fig. \ref{fig:bad_waves_2h}). Contrary, the collocated grid scheme Ah63 fails to adequately reproduce stationary state: the exact smooth feature is modulated by approximately $2\Delta x$ wave. This is due to the (near) zero phase velocity of shortest resoled scale waves that makes the effective Rossby deformation radius for these waves very small. The staggered grid schemes are free of this deficiency.

\section{Conclusion}\label{sec:conclusion}
We presented staggered grid Summation-By-Parts Finite-Difference (SBP FD) provable-stable approach for numerical solution of the wave-dominated problems in closed domains discretized with a multi-block curvilinear grid. This article extends the pioneering works \cite{OReilly_staggered_sbp_2020} and \cite{Gao_staggered_sbp_2019} on staggered SBP FD methods. Unlike \cite{OReilly_staggered_sbp_2020} the presented approach uses a uniform grid without additional collocated nodes at the blocks interfaces. In this aspect our approach is similar to \cite{Gao_staggered_sbp_2019}, but extends it to the case of curvilinear coordinate mappings inside grid blocks and provide higher order accurate operators.

To impose the interface conditions at the grid block edges, we derived the new hybrid Simultaneous-Approximation-Terms (SAT)-Projection method. As compared to the commonly-utilized pure SAT method, the discrete Laplace operator spectrum of SAT-Projection method is free of spurious fast and stationary modes associated with the discontinuity of solution at grid block interfaces. This allows to use significantly larger time-steps (up to twice larger for 6/3 order accurate scheme) and avoid accumulation of grid-scale noise near interfaces.

We also introduce energy-conserving Coriolis terms approximations, extending the applicability of staggered SBP-FD method to geophysical hydrodynamics problems. The technique to ensure continuity of Coriolis terms across grid blocks interfaces is also presented.

We test the presented staggered SBP-FD numerical method using linearized shallow water equations on the cubed-sphere grid. This approximated model of atmospheric and oceanic dynamics is very indicative for testing how numerical methods are representing such basic dynamics as inertia-gravity waves propagation, Rossby-waves dynamics and Rossby-adjustment that lies behind the more complex geophysical hydrodynamics phenomena.

One can see that staggered SBP-FD schemes outperform collocated SBP-FD in terms of reproduction accuracy for pure gravity waves of various spatial scales. Staggered schemes do not significantly retard the waves as short as $4\Delta x$ that is not the case for collocated grid schemes.

The presented energy-conserving Coriolis terms approximations involve the interpolation operators between the grid locations of velocity components. This is an additional source of errors as compared to collocated grid schemes where Coriolis terms are trivial. Despite this, the presented staggered grid schemes are at least as good as collocated counterparts in the reproduction of large-scale inertia-gravity waves and Rossby-adjusted features. For noisy initial fields, staggered schemes outperforms the collocated ones in reproduction of Rossby adjustment process as expected from \cite{Randall_GeostAdj}.

Despite its advantages, the presented method has certain limitations. 
As compared to more common collocated SBP-schemes, the payoff for better accuracy and dispersion is the computational overhead associated with the interpolations between staggered grid fields locations. 
Also, unlike unstaggered case, stability of the method depends on the skewness of a curvilinear grid.
Additionally, while the method was shown to be effective in handling geophysical waves in closed domain without boundaries, further research is needed to optimize its performance for boundary conditions, such as those found in ocean equations.

To conclude, the numerical solutions to the linearized shallow water equations from the presented staggered SBP-FD schemes are stable, high-order accurate and show good dispersion properties for short-scale waves. Therefore, the experiments with staggered SBP-FD method confirm the properties one can expect from its theoretical basis.

Future research could focus on extending the presented staggered SBP FD framework to more complex systems, such as fully nonlinear shallow water equations or 3D primitive hydrostatic and non-hydrostatic systems relevant for atmospheric models. 

\section*{Data availability}
The code of shallow water model is available at \url{https://github.com/vvshashkin/ch_swe}.

\bibliographystyle{unsrt}

\begin{thebibliography}{10}
	
	\bibitem{kreiss1974finite}
	H-O Kreiss and Godela Scherer.
	\newblock Finite element and finite difference methods for hyperbolic partial
	differential equations.
	\newblock In {\em Mathematical aspects of finite elements in partial
		differential equations}, pages 195--212. Elsevier, 1974.
	
	\bibitem{strand_sbp_1994}
	Bo~Strand.
	\newblock Summation by parts for finite difference approximations for d/dx.
	\newblock {\em Journal of Computational Physics}, 110(1):47--67, 1994.
	
	\bibitem{Olsson_projection_I}
	Pelle Olsson.
	\newblock Summation by parts, projections, and stability {I}.
	\newblock {\em Mathematics of Computation}, 64(211):1035--1065, 1995.
	
	\bibitem{Olsson_projection_II}
	Pelle Olsson.
	\newblock Summation by parts, projections, and stability. {II}.
	\newblock {\em Mathematics of Computation}, 64(212):1473--1493, 1995.
	
	\bibitem{SBP_rev_2014}
	David~C. {Del Rey Fern\'andez}, Jason~E. Hicken, and David~W. Zingg.
	\newblock Review of summation-by-parts operators with simultaneous
	approximation terms for the numerical solution of partial differential
	equations.
	\newblock {\em Computers {\&} Fluids}, 95:171--196, 2014.
	
	\bibitem{lundquist2024encapsulated}
	Tomas Lundquist, Andrew~R Winters, and Jan Nordstr{\"o}m.
	\newblock Encapsulated generalized summation-by-parts formulations for
	curvilinear and non-conforming meshes.
	\newblock {\em Journal of Computational Physics}, 498:112699, 2024.
	
	\bibitem{svard2014review}
	Magnus Sv{\"a}rd and Jan Nordstr{\"o}m.
	\newblock Review of summation-by-parts schemes for initial--boundary-value
	problems.
	\newblock {\em Journal of Computational Physics}, 268:17--38, 2014.
	
	\bibitem{svard2004coordinate}
	Magnus Sv{\"a}rd.
	\newblock On coordinate transformations for summation-by-parts operators.
	\newblock {\em Journal of Scientific Computing}, 20:29--42, 2004.
	
	\bibitem{gustafsson1975convergence}
	Bertil Gustafsson.
	\newblock The convergence rate for difference approximations to mixed initial
	boundary value problems.
	\newblock {\em Mathematics of Computation}, 29(130):396--406, 1975.
	
	\bibitem{mattsson2010stableMultiblock}
	Ken Mattsson and Mark~H Carpenter.
	\newblock Stable and accurate interpolation operators for high-order multiblock
	finite difference methods.
	\newblock {\em SIAM Journal on Scientific Computing}, 32(4):2298--2320, 2010.
	
	\bibitem{almquist2019order}
	Martin Almquist, Siyang Wang, and Jonatan Werpers.
	\newblock Order-preserving interpolation for summation-by-parts operators at
	nonconforming grid interfaces.
	\newblock {\em SIAM Journal on Scientific Computing}, 41(2):A1201--A1227, 2019.
	
	\bibitem{tretyak2023multiresolution}
	Ilya~D Tretyak, Gordey~S Goyman, and Vladimir~V Shashkin.
	\newblock Multiresolution approximation for shallow water equations using
	summation-by-parts finite differences.
	\newblock {\em Russian Journal of Numerical Analysis and Mathematical
		Modelling}, 38(6):393--407, 2023.
	
	\bibitem{carpenter1994time}
	Mark~H Carpenter, David Gottlieb, and Saul Abarbanel.
	\newblock Time-stable boundary conditions for finite-difference schemes solving
	hyperbolic systems: methodology and application to high-order compact
	schemes.
	\newblock {\em Journal of Computational Physics}, 111(2):220--236, 1994.
	
	\bibitem{mattsson2018improvedProjection}
	Ken Mattsson and Pelle Olsson.
	\newblock An improved projection method.
	\newblock {\em Journal of Computational Physics}, 372:349--372, 2018.
	
	\bibitem{ArakawaLamb1977}
	A~Arakawa and V.~Lamb.
	\newblock {\em Computational design of the basic dynamical processes of the
		UCLA general circulation model}, volume~17, pages 173--265.
	\newblock New York: Academic Press, 1977.
	
	\bibitem{Randall_GeostAdj}
	D.~Randall.
	\newblock {Geostrophic ~adjustment ~and ~the ~finite-difference ~shallow ~water
		~equations}.
	\newblock {\em Mon. Wea. Rev.}, 122:1371--1377, 1994.
	
	\bibitem{Thuburn2011_hor_basic_ideas}
	John Thuburn.
	\newblock {\em Horizontal Discretizations: Some Basic Ideas}, pages 43--57.
	\newblock Springer Berlin Heidelberg, Berlin, Heidelberg, 2011.
	
	\bibitem{zingg2000comparison}
	David~W Zingg.
	\newblock Comparison of high-accuracy finite-difference methods for linear wave
	propagation.
	\newblock {\em SIAM Journal on Scientific Computing}, 22(2):476--502, 2000.
	
	\bibitem{GoymanShashkinJCP2021}
	Gordey~S Goyman and Vladimir~V Shashkin.
	\newblock Horizontal approximation schemes for the staggered reduced
	latitude-longitude grid.
	\newblock {\em Journal of Computational Physics}, 434:110234, 2021.
	
	\bibitem{OReilly_staggered_sbp_2017}
	Ossian O'Reilly, Tomas Lundquist, Eric~M. Dunham, and Jan Nordström.
	\newblock Energy stable and high-order-accurate finite difference methods on
	staggered grids.
	\newblock {\em Journal of Computational Physics}, 346:572--589, 2017.
	
	\bibitem{OReilly_staggered_sbp_2020}
	Ossian O'Reilly and N.~Anders Petersson.
	\newblock Energy conservative {SBP} discretizations of the acoustic wave
	equation in covariant form on staggered curvilinear grids.
	\newblock {\em Journal of Computational Physics}, 411:109386, 2020.
	
	\bibitem{walters1983analysis}
	Roy~A Walters and Graham~F Carey.
	\newblock Analysis of spurious oscillation modes for the shallow water and
	{N}avier-{S}tokes equations.
	\newblock {\em Computers \& Fluids}, 11(1):51--68, 1983.
	
	\bibitem{le2005some}
	DY~Le~Roux, A~Sene, V~Rostand, and Emmanuel Hanert.
	\newblock On some spurious mode issues in shallow-water models using a linear
	algebra approach.
	\newblock {\em Ocean Modelling}, 10(1-2):83--94, 2005.
	
	\bibitem{Stan_grid_rev}
	A.~Staniforth and J.~Thuburn.
	\newblock Horizontal grids for global weather and climate prediction models: a
	review.
	\newblock {\em Quart. J. Roy. Met. Soc.}, 138:1 -- 26, 2012.
	
	\bibitem{Gao_staggered_sbp_2019}
	Longfei Gao, David~C. {Del Rey Fern\'andez}, Mark Carpenter, and David Keyes.
	\newblock {SBP-SAT} finite difference discretization of acoustic wave
	equations on staggered block-wise uniform grids.
	\newblock {\em Journal of Computational and Applied Mathematics}, 348:421--444,
	2019.
	
	\bibitem{Williamson_tests}
	D.L. Williamson, J.B. Drake, J.J. Hack, R.~Jakob, and P.N. Swartztrauber.
	\newblock A standard test set for numerical approximations to the shallow water
	equations in spherical geometry.
	\newblock {\em J. Comput. Phys.}, 102:211 -- 224, 1992.
	
	\bibitem{Galewsky_test}
	J.~Galewsky, R.K. Scott, and L.M. Polvani.
	\newblock An initial value problem for testing numerical models of the global
	shallow water equations.
	\newblock {\em Tellus A}, 56:429 -- 440, 2004.
	
	\bibitem{Diener_SBP_cubsph_like}
	Peter Diener, Ernst Dorband, Erik Schnetter, and Manuel Tiglio.
	\newblock Optimized high-order derivative and dissipation operators satisfying
	summation by parts, and applications in three-dimensional multi-block
	evolutions.
	\newblock {\em J. Sci. Comput.}, 32:109--145, 06 2007.
	
	\bibitem{eriksson2023boundary}
	Gustav Eriksson, Jonatan Werpers, David Niemel{\"a}, Niklas Wik, Valter
	Zethrin, and Ken Mattsson.
	\newblock Boundary and interface methods for energy stable finite difference
	discretizations of the dynamic beam equation.
	\newblock {\em Journal of Computational Physics}, 476:111907, 2023.
	
	\bibitem{eriksson2023non}
	Gustav Eriksson.
	\newblock Non-conforming interface conditions for the second-order wave
	equation.
	\newblock {\em Journal of Scientific Computing}, 95(3):92, 2023.
	
	\bibitem{Rancic_CSSW_1996}
	M.~Rančić, R.~J. Purser, and F.~Mesinger.
	\newblock A global shallow-water model using an expanded spherical cube:
	Gnomonic versus conformal coordinates.
	\newblock {\em Quarterly Journal of the Royal Meteorological Society},
	122(532):959--982, 1996.
	
	\bibitem{ThuburnWoolings2005}
	John Thuburn and T.J. Woollings.
	\newblock Vertical discretizations for compressible {E}uler equation
	atmospheric models giving optimal representation of normal modes.
	\newblock {\em Journal of Computational Physics}, 203:386--404, 03 2005.
	
	\bibitem{Ullrich_SW_2010}
	Paul~A. Ullrich, Christiane Jablonowski, and Bram {van Leer}.
	\newblock High-order finite-volume methods for the shallow-water equations on
	the sphere.
	\newblock {\em Journal of Computational Physics}, 229(17):6104 -- 6134, 2010.
	
	\bibitem{ShashkinGoymanExpSL}
	Vladimir~V. Shashkin and Gordey~S. Goyman.
	\newblock Semi-{L}agrangian exponential time-integration method for the shallow
	water equations on the cubed sphere grid.
	\newblock {\em Russian Journal of Numerical Analysis and Mathematical
		Modelling}, 35(6):355--366, 2020.
	
	\bibitem{Melvin_mixedfe_2019}
	Thomas Melvin, Tommaso Benacchio, Ben Shipway, Nigel Wood, John Thuburn, and
	Colin Cotter.
	\newblock A mixed finite-element, finite-volume, semi-implicit discretization
	for atmospheric dynamics: Cartesian geometry.
	\newblock {\em Quarterly Journal of the Royal Meteorological Society},
	145(724):2835--2853, 2019.
	
	\bibitem{thuburn2007rossby}
	J~Thuburn.
	\newblock Rossby wave dispersion on the {C}-grid.
	\newblock {\em Atmospheric Science Letters}, 8(2):37--42, 2007.
	
	\bibitem{thuburn2004CgridCons}
	John Thuburn and Andrew Staniforth.
	\newblock Conservation and linear rossby-mode dispersion on the spherical {C}
	grid.
	\newblock {\em Monthly Weather Review}, 132(2):641 -- 653, 2004.
	
	\bibitem{horn2012matrix}
	Roger~A Horn and Charles~R Johnson.
	\newblock {\em Matrix analysis}.
	\newblock Cambridge university press, 2012.
	
	\bibitem{ShashkinGoymanTolstykh_SBP_SWE_2022}
	Vladimir~V. Shashkin, Gordey~S. Goyman, and Mikhail~A. Tolstykh.
	\newblock Summation-by-parts finite-difference shallow water model on the
	cubed-sphere grid. part {I}: Non-staggered grid.
	\newblock {\em Journal of Computational Physics}, 474:111797, 2023.
	
\end{thebibliography}

\section*{SUPPLEMENTARY MATERIALS}

\appendix\label{sec:appendix}

\section{Coefficients of SBP matrices}\label{app:coefficients}

Below, we present the coefficients of the staggered summation-by-parts quadrature, differentiation and interpolation matrices used in this work. 

\subsection{Quadrature matrices}
Below, coefficients of quadrature matrices $H_c$ and $H_v$ are presented.

\subsubsection{Second order}

\begin{equation}
    \begin{gathered}
        H_v = \Delta x \cdot \mathrm{diag}\left(\frac{1}{2},1,...,1,\frac{1}{2}\right), \\
        H_c = \Delta x \cdot \mathrm{diag}(1,1,...,1,1).
    \end{gathered}
\end{equation}

\subsubsection{Fourth order}\label{app:coefficients:quad:42}
\begin{equation}
\begin{gathered}
    H_v = \Delta x \cdot \mathrm{diag}\left(\frac{7}{18},\frac{9}{8},1,\frac{71}{72},1,...,1,\frac{71}{72},1,\frac{9}{8},\frac{7}{18}\right), \\
    H_c = \Delta x \cdot \mathrm{diag}\left(\frac{13}{12}, \frac{7}{8}, \frac{25}{24}, 1,...,1,\frac{25}{24},\frac{7}{8},\frac{13}{12}\right).
\end{gathered}
\end{equation}

\subsubsection{Sixth order}\label{app:coefficients:quad:63}

\begin{equation}
\begin{gathered}
    H_v = \Delta x \cdot \mathrm{diag}(h_{1}^v, h_{2}^v, h_{3}^v, h_{4}^v, h_{5}^v,1,...,1,h_{5}^v,h_{4}^v,h_{3}^v,h_{2}^v,h_{1}^v), \\
    h_{1}^v = \frac{95}{288},\ h_{2}^v = \frac{317}{240},\ h_{3}^v = \frac{23}{30},\ h_{4}^v = \frac{793}{720},\ h_{5}^v = \frac{157}{160},\\
    H_c = \Delta x \cdot \mathrm{diag}(h_{1}^c, h_{2}^c, h_{3}^c, h_{4}^c, h_{5}^c,  h_{6}^c,1,...,1,h_{6}^c,h_{5}^c,h_{4}^c,h_{3}^c,h_{2}^c,h_{1}^c), \\
    h_{1}^c = \frac{325363}{276480},\ h_{2}^c = \frac{144001}{276480},\ h_{3}^c = \frac{43195}{27648},\ h_{4}^c = \frac{86857}{138240},\ h_{5}^c = \frac{312623}{276480},\ h_{6}^c = \frac{271229}{276480}.
    \end{gathered}
\end{equation}

\subsection{Differentiation matrices}\label{app:coefficients:diff}
Below, the $D_{cv}$ differentiation matrices and boundary extrapolation operators $\mathbf{r}$, $\mathbf{l}$ are presented. $D_{vc}$ can be obtained using the relation $H_v D_{cv} + D_{vc}^T H_c = \mathbf{e_r}^T\mathbf{r} - \mathbf{e_l}^T\mathbf{l}$, where $\mathbf{e}_r = (1,0,...,0)^T$, $\mathbf{e}_l = (0,...,0,1)^T$.

\subsubsection{Fourth order}\label{app:coefficients:diff:42}
\begin{equation}
		\mathrm D_{cv} = \frac{1}{\Delta x}
\begin{bmatrix}
	\text{-}2 & 3  & \text{-}1 &  0 & 0 & & & \\
         \text{-}1  & 1   &  0  & 0 & 0 & & & \\
	\frac{1}{24} & \text{-}\frac{9}{8} & \frac{9}{8} & \text{-}\frac{1}{24} &0 & & &\\
 \text{-}\frac{1}{71} & \frac{6}{71} & \text{-}\frac{83}{71} & \frac{81}{71} &\text{-}\frac{3}{71} & & &\\
                        &               & \frac{1}{24} & \text{-}\frac{9}{8} & \frac{9}{8} & \text{-}\frac{1}{24} & & \\
                        &               &               & \ddots & \ddots & \ddots & \ddots & 
\end{bmatrix},
\end{equation}
The lower-right corner elements are determined as $d_{(N+2-i) (N+1-j)} = -d_{ij}$.
\begin{equation}
\mathbf{l} =\frac{1}{8}(15, -10, 3,0, \dots, 0)^T 	
	, \quad
		\mathbf{r} =\frac{1}{8}(0, \dots,0, 3, -10, 15)^T.
\end{equation}

\subsubsection{Sixth order}\label{app:coefficients:diff:63}
\begin{equation}
\setcounter{MaxMatrixCols}{20}  
D_{cv} = \frac{1}{\Delta x}
\begin{bmatrix}
          d_{11} & d_{12}  & d_{13} &  d_{14} & d_{15} & d_{16} & & & & &\\
          d_{21} & d_{22}  & d_{23} &  d_{24} & d_{25} & d_{26} & & & & &\\
	   d_{31} & d_{32}  & d_{33} &  d_{34} & d_{35} & d_{36} & & & & &\\
          d_{41} & d_{42}  & d_{43} &  d_{44} & d_{45} & d_{46} & & & & &\\
          d_{51} & d_{52}  & d_{53} &  d_{54} & d_{55} & d_{56} & d_{57} & & & &\\
               0 &      0  & \text{-}\frac{3}{640} & \frac{25}{384} & \text{-}\frac{75}{256} & \frac{75}{256} & \text{-}\frac{25}{384}& \frac{3}{640} &  & & \\
                 &         &        & \ddots & \ddots & \ddots & \ddots & \ddots & \ddots &  & 
\end{bmatrix},
\end{equation}
The lower-right corner elements are determined as $d_{(N+2-i) (N+1-j)} = -d_{ij}$.

\begin{table}[H]
    \centering
    \def\arraystretch{1.3}
    \begin{tabular}{ll}
        $d_{11} = \frac{-60711983+15005904\,{\it c55}+5183400\,{\it c34}}{21888000}$, & $d_{12} = \frac{101173243-30011808\,{\it c55}-15550200\,{\it c34}}{17510400}$,\\
    $d_{13} = \frac{-7780959+1727800{\it c34}}{1459200}$,&  $d_{14} = \frac{-5183400{\it c34}+35609465+30011808{\it c55}}{8755200}$,\\
    $d_{15} = \frac{-7502952{\it c55}-5209847}{2188800}$,& $d_{16} = \frac{18712829+30011808{\it c55}+1727800{\it c34}}{29184000}$, \\
    $d_{21}=\frac{-53376169-30011808\,{\it c55}-10366800\,{\it c34}}{43822080\,}$,&
    $d_{22} = {\frac{7190801+10003936\,{\it c55}+5183400\,{\it c34}}{5842944\,}}$, \\
    $d_{23}=-{\frac{2591700\,{\it c34}-2846555}{2191104\,}}$,&
    $d_{24} = {\frac{-27181195-30011808\,{\it c55}+5183400\,{\it c34}}{8764416\,}}$, \\ 
    $d_{25} = {\frac{7223559+10003936\,{\it c55}}{2921472\,}}$, & 
    $d_{26} = \frac{-59866697-90035424\,{\it c55}-5183400\,{\it c34}}{87644160\,}$,\\
    $d_{31} = {\frac {332488+625246\,{\it c55}+215975\,{\it c34}}{353280\,}}$,&
    $d_{32} = {\frac {-6326795-10003936\,{\it c55}-5183400\,{\it c34}}{2260992\,}}$,\\
    $d_{33} = {\frac {1727800\,{\it c34}-665205}{565248\,}}$,&
    $d_{34} = {\frac {8940511+10003936\,{\it c55}-1727800\,{\it c34}}{1130496\,}}$,\\
    $d_{35} = {\frac {-3843253-5001968\,{\it c55}}{565248\,}}$,&
    $d_{36} = {\frac {21758409+30011808\,{\it c55}+1727800\,{\it c34}}{11304960\,}}$,\\
    $d_{41} = {\frac {-17586239-30011808\,{\it c55}-10366800\,{\it c34}}{36541440\,}}$,&
    $d_{42} ={\frac {14084351+30011808\,{\it c55}+15550200\,{\it c34}}{14616576\,}}$,\\
$d_{43} =-{\frac {215975\,{\it c34}}{152256\,}}$,&
$d_{44} ={\frac {5183400\,{\it c34}-30011808\,{\it c55}-21697151}{7308288\,}}$,\\
$d_{45} ={\frac {25503551+30011808\,{\it c55}}{7308288\,}}$,&
$d_{46} ={\frac {-24437759-30011808\,{\it c55}-1727800\,{\it c34}}{24360960\,}}$, \\
$d_{51} = {\frac {9606527+15005904\,{\it c55}+5183400\,{\it c34}}{65111040\,}}$,&
$d_{52} = {\frac {-4598783-10003936\,{\it c55}-5183400\,{\it c34}}{17362944\,}}$,\\
$d_{53} = {\frac {-4811905+5183400\,{\it c34}}{13022208\,}}$,&
$d_{54} = {\frac{5665537-5183400\,{\it c34}+30011808\,{\it c55}}{26044416\,}}$,\\
$d_{55} = -{\frac {312623\,{\it c55}}{271296\,}}$,&
$d_{56} = {\frac {68894207+90035424\,{\it c55}+5183400\,{\it c34}}{260444160\,}}$,\\
$d_{57} = {\frac {3}{628\,}}$
    \end{tabular}
\end{table}

Here, $c_{34}$ and $c_{55}$ are free parameters. In this work, we use two sets of parameters. The first set 
\begin{align}
    {\it c34} =0.6690374220138081 , \ {\it c14} = -0.7930390145751754 ,
\end{align}
is obtained by the minimization of \eqref{eq:poly_opt_63} objective function. While the second one
\begin{align}
    {\it c34} = 0.467391226104632 , \ {\it c55} = -0.723617281756727 ,
\end{align}
is obtained by the minimization of \eqref{eq:wave_opt_63} objective function.

\begin{equation}
    \mathbf{l} = \frac{1}{16}(35, -35, 21, -5,0,...,0)^T, \quad \mathbf{r} = \frac{1}{16}(0,...,0,-5,21,-35,35)^T.
\end{equation}

\subsection{Interpolation matrices}\label{app:coefficients:interp}

Below, the $P_{vc}$ interpolation matrices are presented. $P_{cv}$ can be obtained using the relation $H_v P_{cv} = P_{vc}^T H_c$.  

\subsubsection{Fourth order}\label{app:coefficients:interp:42}

\begin{equation}
		P_{vc} =
\begin{bmatrix}
          p_{11} & p_{12}  & p_{13} &  p_{14} & 0 & & & \\
          p_{21} & p_{22}  & p_{23} &  p_{24} & 0 & & & \\
	   p_{31} & p_{32}  & p_{33} &  p_{34} & p_{35} & & &\\
               0 &      0  & \frac{\text{-}1}{16} & \frac{9}{16} & \frac{9}{16} & \frac{\text{-}1}{16} & & \\
                 &         &        & \ddots & \ddots & \ddots & \ddots & 
\end{bmatrix},
\end{equation}
The lower-right corner elements are determined as $p_{(N+1-i) (N+2-j)} = p_{ij}$.

The coefficients in $P_{vc}$ are given by:
\begin{table}[H]
    \centering
    \begin{tabular}{ll}
        $p_{11} = \frac{1}{2}+{\it c13}+2\,{\it c14}$ & $p_{12} = \frac{1}{2}-2\,{\it c13}-3\,{\it c14}$ \\
        $p_{13} = {\it c13}$ & $p_{14} = {\it c14}$ \\
        $p_{21} = -{\frac{8}{63}}-{\frac {52\,{\it c13}}{21}}-{\frac {104\,{\it c14}}{21}}$ & $p_{22} = {\frac{29}{42}}+{\frac {104\,{\it c13}}{21}}+{\frac {52\,{\it c14}}{7}}$ \\
        $p_{23} = -{\frac {52\,{\it c13}}{21}}+\frac{1}{2}$ & $p_{24} = -{\frac{4}{63}}-{\frac {52\,{\it c14}}{21}}$ \\
        $p_{31} = {\frac {26\,{\it c13}}{25}}+{\frac {52\,{\it c14}}{25}}-\frac{1}{25}$ & $p_{32} = -{\frac{1}{50}}-{\frac {52\,{\it c13}}{25}}-{\frac {78\,{\it c14}}{25}}$ \\
        $p_{33} = \frac{3}{5}+{\frac {26\,{\it c13}}{25}}$ & $p_{34} = {\frac{13}{25}}+{\frac {26\,{\it c14}}{25}}$ \\
        $p_{35} = -{\frac{3}{50}}$ &   \\
    \end{tabular}
\end{table}
Here, ${\it c13}$ and ${\it c14}$ are free parameters. In this work, we use the following values
\begin{align}
    {\it c13} =\frac{102207746025903}{808013506696916} , \ {\it c14} = -\frac{289843969221617}{9696162080362992} ,
\end{align}
obtained by the minimization of (\ref{eq:interp_42_obj_func}) objective function.

\subsubsection{Sixth order}\label{app:coefficients:interp:63}

\begin{equation}
\setcounter{MaxMatrixCols}{20}  
P_{vc} = 
\begin{bmatrix}
          p_{11} & p_{12}  & p_{13} &  p_{14} & p_{15} & & & & & &\\
          p_{21} & p_{22}  & p_{23} &  p_{24} & p_{25} & & & & & &\\
	   p_{31} & p_{32}  & p_{33} &  p_{34} & p_{35} & p_{36} & & & & &\\
          p_{41} & p_{42}  & p_{43} &  p_{44} & p_{45} & p_{46} & p_{47} & & & &\\
          p_{51} & p_{52}  & p_{53} &  p_{54} & p_{55} & p_{56} & p_{57} & p_{58} & & &\\
          p_{61} & p_{62}  & p_{63} &  p_{64} & p_{65} & p_{66} & p_{67} & p_{68} & p_{69} & &\\
               0 &      0  & 0 & 0 & \frac{3}{256} & \frac{\text{-}25}{256} & \frac{150}{256}& \frac{150}{256} & \frac{\text{-}25}{256} & \frac{3}{256} & \\
                 &         &        &  & & \ddots & \ddots & \ddots & \ddots & \ddots & \ddots 
\end{bmatrix},
\end{equation}
The lower-right corner elements are determined as $p_{(N+1-i) (N+2-j)} = p_{ij}$.

The coefficients in $P_{vc}$ are given by:
\begin{table}[H]
    \centering
    \def\arraystretch{1.3}
    \begin{tabular}{l}
        $p_{11} = {\frac{4474753}{7808712}}+{\frac {312623\,{\it c53}}{650726}}+{\frac {937869\,{\it c52}}{650726}}-{\frac {813687\,{\it c64}}{1301452}}+{\frac {2441061\,{\it c62}}{1301452}}+{\frac {86857\,{\it c42}}{325363}}+{\frac {86857\,{\it c43}}{976089}}$, \\
        $p_{12} = {\frac{136944}{325363}}-{\frac {173714\,{\it c42}}{325363}}-{\frac {1627374\,{\it c62}}{325363}}-{\frac {937869\,{\it c52}}{325363}}$, \\
        $p_{13} = -{\frac{848457}{2602904}}-{\frac {173714\,{\it c43}}{325363}}+{\frac {2441061\,{\it c64}}{650726}}+{\frac {2441061\,{\it c62}}{650726}}-{\frac {937869\,{\it c53}}{325363}}$, \\
        $p_{14} = {\frac{2331127}{3904356}}+{\frac {173714\,{\it c42}}{325363}}+{\frac {694856\,{\it c43}}{976089}}-{\frac {1627374\,{\it c64}}{325363}}+{\frac {1250492\,{\it c53}}{325363}}+{\frac {937869\,{\it c52}}{325363}}$, \\
        $p_{15} = -{\frac{10145}{38278}}-{\frac {937869\,{\it c53}}{650726}}-{\frac {937869\,{\it c52}}{650726}}+{\frac {2441061\,{\it c64}}{1301452}}-{\frac {813687\,{\it c62}}{1301452}}-{\frac {86857\,{\it c42}}{325363}}-{\frac {86857\,{\it c43}}{325363}}$, \\
        $p_{21} = -{\frac{373145}{354464}}-{\frac {86857\,{\it c43}}{144001}}+{\frac {4068435\,{\it c64}}{1152008}}-{\frac {12205305\,{\it c62}}{1152008}}-{\frac {1250492\,{\it c53}}{432003}}-{\frac {260571\,{\it c42}}{144001}}-{\frac {1250492\,{\it c52}}{144001}}$, \\
        $p_{22} ={\frac{273888}{144001}}+{\frac {521142\,{\it c42}}{144001}}+{\frac {4068435\,{\it c62}}{144001}}+{\frac {2500984\,{\it c52}}{144001}}$, \\
        $p_{23} ={\frac{520551}{209456}}+{\frac {521142\,{\it c43}}{144001}}-{\frac {12205305\,{\it c64}}{576004}}-{\frac {12205305\,{\it c62}}{576004}}+{\frac {2500984\,{\it c53}}{144001}}$, \\
        $p_{24} =-{\frac{1142117}{288002}}-{\frac {521142\,{\it c42}}{144001}}-{\frac {694856\,{\it c43}}{144001}}+{\frac {4068435\,{\it c64}}{144001}}-{\frac {10003936\,{\it c53}}{432003}}-{\frac {2500984\,{\it c52}}{144001}}$, \\
        $p_{25} ={\frac{7516251}{4608032}}+{\frac {260571\,{\it c43}}{144001}}-{\frac {12205305\,{\it c64}}{1152008}}+{\frac {4068435\,{\it c62}}{1152008}}+{\frac {1250492\,{\it c53}}{144001}}+{\frac {260571\,{\it c42}}{144001}}+{\frac {1250492\,{\it c52}}{144001}}$, \\
        $p_{31} = {\frac{930131}{6911200}}+{\frac {312623\,{\it c53}}{431950}}-{\frac {271229\,{\it c64}}{345560}}+{\frac {86857\,{\it c43}}{431950}}+{\frac {813687\,{\it c62}}{345560}}+{\frac {937869\,{\it c52}}{431950}}+{\frac {260571\,{\it c42}}{431950}}$,\\
        $p_{32} = -{\frac{22824}{215975}}-{\frac {260571\,{\it c42}}{215975}}-{\frac {271229\,{\it c62}}{43195}}-{\frac {937869\,{\it c52}}{215975}}$, \\ 
        $p_{33} = -{\frac{424911}{3455600}}-{\frac {937869\,{\it c53}}{215975}}+{\frac {813687\,{\it c64}}{172780}}-{\frac {260571\,{\it c43}}{215975}}+{\frac {813687\,{\it c62}}{172780}}$, \\
        $p_{34} = {\frac{330902}{215975}}+{\frac {1250492\,{\it c53}}{215975}}+{\frac {937869\,{\it c52}}{215975}}-{\frac {271229\,{\it c64}}{43195}}+{\frac {260571\,{\it c42}}{215975}}+{\frac {347428\,{\it c43}}{215975}}$, \\
        $p_{35} = -{\frac{615889}{1382240}}-{\frac {937869\,{\it c53}}{431950}}+{\frac {813687\,{\it c64}}{345560}}-{\frac {260571\,{\it c43}}{431950}}-{\frac {271229\,{\it c62}}{345560}}-{\frac {937869\,{\it c52}}{431950}}-{\frac {260571\,{\it c42}}{431950}}$, \\
    \end{tabular}
    \centering
    \begin{tabular}{lll}
    $p_{36} = {\frac{324}{43195}}$, &
$p_{41} = -{\frac{17737}{4169136}}-\frac{3{\it c42}+{\it c43}}{6}$, &
$p_{42} = {\it c42}$, \\
$p_{43} = {\it c43}$, &
$p_{44} = {\frac{415759}{1042284}}-\frac{3{\it c42} + 4{\it c43}}{3}$, &
$p_{45} = {\frac{1031359}{1389712}}+\frac{{\it c42}}{2}+\frac{{\it c43}}{2}$, \\
$p_{46} = -{\frac{13500}{86857}}$, &
$p_{47} = {\frac{1620}{86857}}$, &
$p_{51} = {\frac{44783}{5001968}}-\frac{{\it c53 + 3{\it c52}}}{6}$, \\
$p_{52} = {\it c52}$, &
$p_{53} = {\it c53}$, &
$p_{54} = -{\frac{199149}{1250492}}-\frac{4{\it c53}+3{\it c52}}{3}$, \\
$p_{55} = {\frac{3541941}{5001968}}+\frac{{\it c53}}{2}+\frac{{\it c52}}{2}$, &
$p_{56} = {\frac{162000}{312623}}$, &
$p_{57} = -{\frac{27000}{312623}}$, \\
$p_{58} = {\frac{3240}{312623}}$, &
$p_{61} = -{\frac{123231}{8679328}}+\frac{{\it c64}}{8}-\frac{3}{8}\,{\it c62}$, &
$p_{62} = {\it c62}$, \\
$p_{63} = {\frac{30309}{619952}}-\frac{3}{4}\,{\it c64}-\frac{3}{4}\,{\it c62}$, &
$p_{64} = {\it c64}$, &
$p_{65} = -{\frac{1229447}{8679328}}-\frac{3}{8}\,{\it c64}+\frac{{\it c62}}{8}$,\\
$p_{66} = {\frac{162000}{271229}}$,&
$p_{67} = {\frac{162000}{271229}}$,&
$p_{68} = -{\frac{27000}{271229}}$,\\
$p_{69} = {\frac{3240}{271229}}$.
    \end{tabular}
\end{table}
Here, ${\it c42}$, ${\it c43}$, ${\it c52}$, ${\it c53}$, ${\it c62}$, ${\it c64}$  are free parameters. In this work, we use the following values
\begin{equation}
\begin{gathered}
    {\it c42} = -0.3332211159670528, \ {\it c43} = 0.3310769312612241, \ {\it c52} = -0.07099703081266314, \\
    {\it c53} = -0.2916164053358880, \ {\it c62} = 0.05753938634775091, \ {\it c64} = -0.1230378129758785
\end{gathered}
\end{equation}
obtained by the successive minimization of \eqref{eq:interp_63_obj_func_1} and \eqref{eq:interp_63_obj_func_2} objective functions.

\end{document}